\documentclass[12pt]{iopart}

  \expandafter\let\csname equation*\endcsname\relax
  \expandafter\let\csname endequation*\endcsname\relax
\usepackage{iopams, mathdots,times,verbatim,theorem, marginnote,bm,mathtools, setspace, graphicx,marginnote}
\usepackage{iopams}
\usepackage{pgfplots}
\usepackage{graphicx,epsfig,hyperref,verbatim,caption}
\usepackage[labelformat=simple]{subcaption}

\newcommand{\diag}{\mbox{{\rm diag}}}
\newcommand{\R}{\mbox{$\mathbb{R}$}}

\newcommand{\T}{\mbox{\bf T}}

\newcommand{\C}{\mbox{$\mathbb{C}$}}
\newcommand{\Z}{\mbox{$\mathbb{Z}$}}

\newcommand{\F}{\mbox{$\mathbb{F}$}}

\newcommand{\N}{\mbox{$\mathbb{N}$}}
\newcommand{\Sone}{\mbox{\bf S}^1}
\newcommand{\D}{\mbox{$\mathbb{D}$}}

\newcommand{\im}{{\rm Im}}

\newcommand{\Fix}{\mbox{{\rm Fix}}}
\newcommand{\ov}{\overline}

\renewcommand{\dfrac}{\displaystyle\frac}
\newcommand{\proof}{\noindent{\bf Proof:} \quad}

\newcommand{\qed}{\hfill\mbox{\raggedright\rule{0.07in}{0.1in}}\vspace{0.1in}}

\newcommand{\Ref}[1]{(\ref{#1})}
\newcommand{\ds}{\displaystyle}

\newcommand{\AND}{\quad\mbox{and}\quad}

\newtheorem{theorem}{Theorem}[section]
\newtheorem{proposition}[theorem]{Proposition}

\theorembodyfont{\rmfamily}

\begin{document}
\title[]{$\D_n$ Symmetric Hamiltonian System: A Network of Coupled Gyroscopes as a Case Study}

\author{Pietro-Luciano Buono$^a$, Bernard S Chan$^b$,  Antonio Palacios$^b$ and Visarath In$^c$}

\address{$^a$Faculty of Science, University of Ontario Institute of Technology \\
 2000 Simcoe St N, Oshawa, ON L1H 7K4, Canada\\
  $^b$Nonlinear Dynamical Systems Group, Department of Mathematics, \\
  San Diego State University, San Diego, CA 92182 \\
  $^c$Space and Naval Warfare Systems Center, Code 2363, \\
  53560 Hull Street, San Diego, CA 92152-5001, USA}
\ead{luciano.buono@uoit.ca}
\begin{abstract}
The evolution of a large class of biological, physical and engineering systems can be studied through both dynamical systems theory and Hamiltonian mechanics. The former theory, in particular its specialization to study systems with symmetry, is already well developed and has been used extensively on a wide variety of spatio-temporal systems. There are, however, fewer results on higher-dimensional Hamiltonian systems with symmetry. This lack of results has lead us to investigate the role of symmetry, in particular dihedral symmetry, on high-dimensional coupled Hamiltonian systems. As a representative example, we consider the model equations of a ring of vibratory gyroscopes. The equations are reformulated in a Hamiltonian structure and the corresponding  normal forms are derived.
Through a normal form analysis, we investigated the effects of various coupling schemes and unraveled the nature of the bifurcations that lead the ring of gyroscopes into and out of synchronization. The Hamiltonian approach is specially useful in investigating the collective behavior of small and large ring sizes and it can be readily extended to other symmetry-related systems.

\end{abstract}

\maketitle

\section{Introduction}
\label{sec:introduction}
High-dimensional nonlinear systems with symmetry arise naturally at various length scales. Examples
can be found in molecular dynamics~\cite{tuckerman1991molecular}, underwater vehicle dynamics~\cite{leonard1997stability}, magnetic- and electric-field sensors~\cite{bulsara2008exploiting,in2003coupling,in2006complex,palacios2006cooperative}, gyroscopic~\cite{acar2009mems,nagata1998bifurcations} and navigational systems~\cite{grewal2007global,rogers2003applied}, hydroelastic rotating systems~\cite{mcdonald1999non,mcdonald2002parametric,mcdonald2006bifurcations,nagata1998bifurcations},
and  {\em complex systems} such as telecommunication infrastructures~\cite{kocarev2005complex} and
power grids~\cite{susuki2009hybrid,susuki2008global}. Whereas the theory of symmetry breaking
bifurcations of typical invariant sets, i.e., equilibria, periodic solutions, and chaos, is well-developed for
general low-dimensional systems~\cite{golubitsky1988groups,strogatz2001nonlinear}, there are
significantly fewer results on the corresponding theory for symmetric high-dimensional nonlinear
{\em mechanical and electrical} systems, including coupled Hamiltonian systems~\cite{hampton1999measure,skokos2001stability,smereka1998synchronization}. Thus, we aim this work
at advancing the study of the role of symmetry in high-dimensional nonlinear systems with Hamiltonian
structure. We consider systems whose symmetries are represented by the dihedral group $\D_N$,
which describes the symmetries of an $N$-gon, as it arises commonly in generic versions of coupled
network systems with bi-directional coupling. The cyclic group $\Z_N$, which describes networks with
nearest-neighbor coupling with a preferred orientation, i.e., unidirectional coupling, is also described.

As a case study, we consider the model equations of a ring of vibratory gyroscopes. Each gyroscope
is modeled by a 4-dimensional nonautonomous system of Ordinary Differential Equations (ODEs).
Then a network of $N$ gyroscopes is governed by a coupled nonautonomous ODE system of
dimension $4N$, which can be difficult to study when $N$ is large. Numerical simulations show that
under certain conditions, which depend mainly on the coupling strength, the dynamics of the individual
gyroscopes will synchronize with one another~\cite{vu2010two}. A two-time scale analysis, carried out
for the particular case of $N=3$ gyroscopes, yielded an approximate analytical expression for a critical
coupling strength at which the gyroscopic oscillations merge in a pitchfork bifurcation; passed this
critical coupling the synchronized state becomes locally asymptotically stable. The synchronization
pattern is of particular interest because it can lead to a reduction in the phase drift that typically affects
the performance of most gyroscopes. For larger arrays, numerical simulations show that there still
exists a critical value of coupling strength that leads to synchronization and, potentially, to additional
reductions in phase drift. Thus finding an approximate expression for that critical coupling is a very
important task. One possible approach to carry out this task is to generalize the two-time scale
analysis to any $N$. The system of partial differential equations that results from this approach is,
however, too cumbersome and not amenable to analysis. Furthermore, one may have to perform
multiple versions of the same analysis in order to distinguish the different types of bifurcations that
may occur for various combinations of $N$ values. An alternative approach is first to cast the
equations of motion, without forcing, in Hamiltonian form and then study whether the coupled ring
system preserves the Hamiltonian structure. If it does, then, in principle, we could calculate a general
Hamiltonian function, valid for any ring size $N$, from which we can readily determine the existence
of equilibria and their spectral properties.
More importantly, it should also be possible to uncover the critical value of coupling strength that leads
a ring of any size to synchronization and to better understand the nature of the bifurcations for larger $N$.
Finally, the existence of synchronous periodic solutions, its stability and bifurcations at a critical coupling
strength
can be investigated by treating the time-dependent forcing term as a small perturbation of the
Hamiltonian structure.

In this manuscript we show that the second approach outlined above, i.e., via Hamiltonian dynamics,
can indeed provide a more rigorous framework to study the collective behavior of the coupled
gyroscope system. In fact, we show that a coupled ring with $\Z_N$-symmetry does not have a
Hamiltonian structure while a $\D_N$ symmetric ring does admit a Hamiltonian structure. In this latter
case, the Hamiltonian analysis provides, through a normal form analysis and the Equivariant Splitting
Lemma~\cite{bridges1993singularity}, a better picture of the nature of the bifurcations for any ring size $N$ and an exact analytical
expression for the critical coupling strength that leads to synchronized behavior, also valid for any $N$.
We wish to emphasize again that the focus of the theoretical work to gyroscopes and to the symmetry
groups $\Z_N$ and $\D_N$ is not exhaustive. Many other high-dimensional coupled Hamiltonian
systems with symmetry can, in principle, be studied through a similar approach. For instance, ongoing
research work on energy harvesting systems, which also attempts to exploit the collective vibrations
of coupled galfenol-based materials to maximize power output, leads to a high-dimensional nonlinear
system whose model equations are very similar to those of the coupled gyroscope system. The gyroscopes and
the energy harvesting system are only two representative examples of high-dimensional coupled
systems that can benefit from a theoretical framework to study Hamiltonian systems with symmetry.

The manuscript is organized as follows. In section~\ref{sec:hamiltonianFormulation}, the fundamental
principles of operation of a gyroscope and their governing  equations are briefly described for
completeness purposes. The governing equations for a 1D ring-array of $N$ linearly coupled
gyroscopes and their Hamiltonian formulation are also introduced. Proofs of the lack of Hamiltonian
structure for a ring with $\Z_N$ symmetry is presented as well as proof of Hamiltonian structure for a
ring with dihedral $\D_N$ symmetry. Section~\ref{sec:linear-origin} investigates the effect of the
symmetry on the linearized equations from which we obtain explicit expressions for the eigenvalues
and determine the distribution of the spectrum for all $N$.
In section~\ref{sec:linearAnalysis}, symplectic matrices are calculated to transform the quadratic
part of the Hamiltonian function, which corresponds to the linear part in the original coordinates, to normal
form. In section~\ref{sec:nonlinearNormalForm} the normal form calculations are extended to include
nonlinear terms. In particular, the minimum set of invariant terms necessary for the nonlinear system
to be written in normal form are determined. Finally, the splitting lemma is employed to separate the
degenerate and non degenerate components of the Hamiltonian function, so that only the most
essential nonlinear terms are left for further analysis of the bifurcations in the coupled system. In
section~\ref{sec:bifurcationParameter}, the results of the general theory for arbitrary $N$ are illustrated
to study a 1D ring with $\D_3$-symmetry. In section~\ref{sec:conclusion}, some concluding remarks are
presented.

\section{Hamiltonian Formulation}
\label{sec:hamiltonianFormulation}

\subsection{Single Vibratory Gyroscope} \label{sec:background}
\begin{figure}[htb]
\begin{center}
   \epsfig{file=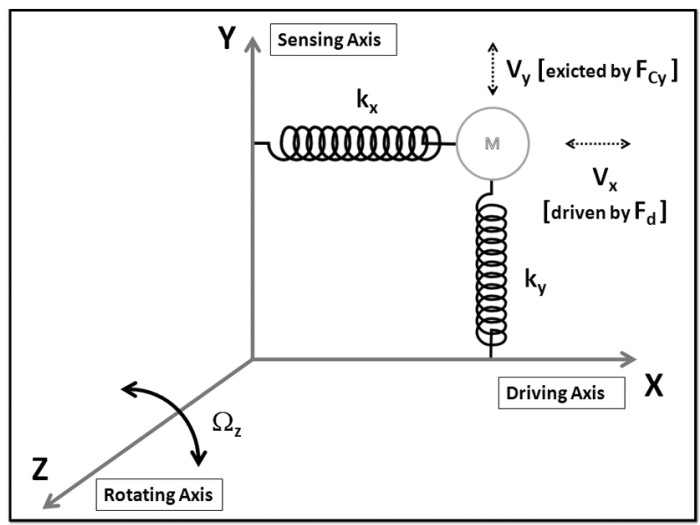,width=4in}
   \caption{Schematic diagram of a vibratory gyroscope system.
   An internal driving force induces the spring-mass system to vibrate in one
   direction, the $x$-axis in this example. An external rotating force, perpendicular
   to the $xy$-plane induces oscillations in the $y$-direction by transferring energy
   through the Coriolis force. These latter oscillations can be used to detect and quantify
   the rate of rotation.}
   \label{fig:mass-spring-model}
\end{center}
\end{figure}

A conventional vibratory gyroscope consists of a proof-mass system as is shown in figure~\ref{fig:mass-spring-model}. The system operates~\cite{acar2009mems,apostolyuk2006theory,apostolyuk2004dynamics}
on the basis of energy transferred from a driving mode to a sensing mode through the Coriolis
force~\cite{CORIOLIS}. In this configuration, a change in the acceleration around the driving $x$-axis
caused by the presence of the Coriolis force induces a vibration in the sensing $y$-axis which can be
converted to measure angular rate output or absolute angles of rotation.
Normally, a higher amplitude response of the $y$-axis translates to  an increase in sensitivity of a gyroscope. Thus, to achieve high sensitivity most gyroscopes operate at resonance in both drive- and sense-modes. But since the phase and frequency of the sense-mode is determined by the phase and frequency of the Coriolis force which itself depends on those of the driving signal, most gyroscopes operate exactly at the drive-mode resonant frequency while the sense-mode frequency is controlled to match the drive-mode resonant frequency. Consequently, the performance of a gyroscope, in terms of accuracy and sensitivity, depends greatly on the ability of the driving signal to produce stable oscillations with constant amplitude, phase, and frequency. To achieve these important requirements, a variety of schemes, based mainly on closed-loops and phase-locked loops circuits, have been proposed~\cite{acar2009mems}. Parametrical resonance in MEMS (Micro-Electro-Mechanical Systems) gyroscopes has also been extensively studied as an alternative to harmonically driven oscillators~\cite{demartini2007linear}. More recently, we have shown that coupling similar gyroscopes in some fashion can lead to globally asymptotically stable synchronized oscillations that are robust enough to mitigate the negative effects of noise while minimizing phase drift~\cite{davies2013collected,vu2010two,vu2011drive}. In those works, perturbation analysis and computer simulations were employed as the main tools to calculate approximate analytical expressions for the boundary curves that separate synchronized behavior from other patterns of collective behavior. However, a better understanding of the underlying bifurcations that lead into and out of synchronization for small and large rings of gyroscopes is still missing. We show in this manuscript that a formulation of the model equations via Hamiltonian functions can provide a better understanding of the role of symmetry in the coupled gyroscope systems and in other generic systems with Hamiltonian structure.

Based on the fundamental principles of operation illustrated in figure~\ref{fig:mass-spring-model},
the governing equations of a single gyroscope can be modeled after a spring-mass system:
\begin{equation} \label{eq:single_gyro}
   \begin{array}{lllll}
        m \ddot{x} + c_x \dot{x} +  \kappa_x x + \mu_x x^3 &=&
        f_e(t) &+& 2 m \Omega_z \dot{y} \\
        m \ddot{y} + c_y \dot{y} + \kappa_y y + \mu_y y^3 &=& &-&
        2 m \Omega_z \dot{x},
  \end{array}
\end{equation}
where $x$ $(y)$ represents the drive (sense) modes, $m$ is the proof mass, $\Omega_z$
is the angular rate of rotation along a perpendicular direction ($z$-axis),
$c_x$ ($c_y$) is the damping coefficient along the $x$- ($y$-) direction, and $\kappa_x$
($\kappa_y$) and $\mu_x$ ($\mu_y$) are the linear and nonlinear damping coefficients
along the $x$- ($y$-) directions, respectively. Typically, the forcing term  has
sinusoidal form $f_e(t)=A_d\cos{w_d t}$. The Coriolis forces appear in the driving- and
sensing-modes as $F_{cx} = 2m\Omega_z \dot{y}$ and $F_{cy} = -2m\Omega_z \dot{x}$,
respectively. Note that the $x$-axis is also excited by a reference driving force where $A_d$
is the amplitude and $w_d$ is the frequency of the excitation. Typical operational parameter
values for this work are shown in table~\ref{tab:parameters}.

\begin{table}[h]
\begin{center}
   {\begin{tabular}{cccc}
    Parameter & Value & Unit            \\  \hline
    $m$                           & 1.0E-09   & Kg          \\[3pt]
    $c_x$,  $c_y$        & 5.1472E-07 & N s/meter \\[3pt]
    $\kappa_x$, $\kappa_y$           & 2.6494     & N/meter \\[3pt]
    $\mu_x$, $\mu_y$ & 2.933       & N/meter$^3$ \\[3pt]
    $A_d$                       & 1.0E-03   & N    \\[3pt]
    $w_d$                      & 5.165E+04 & rad/sec \\
   \end{tabular}}

  \caption{System parameters for a vibratory gyroscope.}
  \label{tab:parameters}
\end{center}
\end{table}

Under these conditions, the gyroscope of figure~\ref{fig:mass-spring-model} can detect an
applied angular rate $\Omega_z$ by measuring the displacements along the $y$-axis
caused by the transfer of energy via the Coriolis force. In the absence of a external rotation,
i.e., $\Omega_z=0$, the equations of motion listed in~\eqref{eq:single_gyro} along the two
axes become uncoupled from one another and the dynamics along the $x$-axis reduces to
that of a Duffing oscillator subject to a periodic force, which has been extensively
studied~\cite{NAYFEH,RAND}. The motion along the $y$-axis, however, eventually
approaches the zero equilibrium due to the damping term $c_y$.

\subsection{The Coupled System}
We now consider $N$ identical gyroscopes coupled identically in a ring system.
We assume that $\mu_x=\mu_y=\mu>0$. Thus, the behavior of the individual $i^{th}$ gyroscope in the
system is described by the following system of differential equations
\begin{equation} \label{eq:coupled_gyro}
   \begin{array}{llll}
      m\Ddot{x}_i+c_x\dot{x}_i + \kappa_x x + \mu_x x_i^3 & = &f_e(t)&+~2m_i \Omega_z \Dot{y}_i+\displaystyle  \sum_{i \sim j} \lambda_{ij}h(x_i,x_j) \\
      m\Ddot{y}_i+c_y\dot{y}_i + \kappa_y y + \mu_y y_i^3  &= & &-~2m_i\Omega_z\Dot{x}_i,
  \end{array}
\end{equation}
where ${i\sim j}$ denotes all the $j^{th}$ gyroscopes that are coupled to the $i^{th}$ gyroscope,
$\lambda_{ij}$ denotes the coupling strength constant, and $h(x_i,x_j)$ is the coupling function,
which depends on the states of the $i^{th}$ gyroscope and all other $j^{th}$ gyroscopes that are
coupled to it. The specific form for $h$ depends on the configuration of the system. In this work,
we consider the case with no damping friction (i.e., $c_x = c_y = 0$) and no \emph{external forcing}, (i.e., $f_e(t)=0$). We suppose the gyroscopes may be linearly coupled to their nearest neighbors. Thus the system can
be configured as a unidirectional or bidirectional ring. The former case leads to a
system of differential equations with $\Z_N$ symmetry, which is the group of cyclic rotations of $N$
objects. The latter case yields a system with $\D_N$  symmetry, which is the group of symmetries
of a regular $N$-gon. The corresponding coupling functions for these two cases are
\[
h(x_{i-1},x_i,x_{i+1}) = x_{i+1}-x_i \AND h(x_{i-1},x_i,x_{i+1}) = (x_{i+1}-x_i) + (x_{i-1}-x_i),
\]
where $i=1,\ldots,N \operatorname{mod} N$. In fact, because the nonlinear terms are given only by
cubic terms each gyroscope is symmetric with respect to $(x_i,y_i)\to (-x_i,-y_i)$.

To write the system in Hamiltonian form, let $q_i = (q_{i1},q_{i2})^T=(x_i,y_i)^T$ be the
configuration components and $p_i = m\dot{q}_i +Gq_i$ be the momentum components, where
\begin{equation*}
G=\left(\begin{array}{cc}
0&-m\Omega
\\
m\Omega&0
\end{array}\right).
\end{equation*}
Directly differentiating the momentum components, we get $ \dot {p}_i = m\ddot{ q}_i + G\dot{ q}_i$. After rearranging terms, we have $m\ddot{ q}_i = \dot{ p}_i - G\dot{ q}_i$. Then the original equations of the
coupled gyroscopic system in~\eqref{eq:coupled_gyro} can be written in the following form
{\small \begin{equation}\label{eq:gen-coupling}
\begin{split}
\left(\begin{array}{c}
\dot{q_i} \\ \dot{p}_i \end{array}\right) ={}&
\left(\begin{array}{cc} -\frac{G}{m} & \frac{1}{m}I_2\\  -(K-\frac{1}{m}G^2 - \lambda \Gamma h_{x_i}(0,0,0)) &  -\frac{G}{m}\end{array}\right)\left(\begin{array}{c} q_i \\ p_i \end{array}\right) \\
{}& + \left(\begin{array}{c} 0 \\ -f_i + \lambda \Gamma \left(h(x_{i-1},x_i,x_{i+1})-h_{x_i}(0,0,0)x_i\right) \end{array}\right),
\end{split}
\end{equation}}
where
\[
\Gamma = \left(\begin{array}{cc} 1 & 0 \\ 0 & 0 \end{array}\right),\quad K = \text{diag}(\kappa_x,\kappa_y), \AND
f_i= \mu \left(\begin{array}{c} x_i^3 \\  y^3_i \end{array}\right).
\]
The symmetry $(x_i,y_i)\to (-x_i,-y_i)$ lifts symplectically to the symmetry $(q_i,p_i)\to (-q_i,-p_i)$
and this symmetry commutes with any permutation. So~(\ref{eq:gen-coupling}) is $\Z_N\times \Z_2$ or $\D_N\times \Z_2$ equivariant.

We now justify the use of the autonomous system~\eqref{eq:gen-coupling} for the existence, stability and bifurcations of the periodically forced system in~\eqref{eq:coupled_gyro}. Let $Z_i = (q_i,p_i )^T$ be the position and momentum components of gyroscope $i$. Then the evolution equations of each individual gyroscope can be represented as
\[
\dot{Z}_i = M_1 Z_i + \left(\begin{array}{c} 0 \\ \lambda \Gamma(q_{i+1}) \end{array}\right) - {F}_i,
\]
where
\[
M_1 = \left(\begin{array}{cc} -\frac{G}{m} & \frac{1}{m}I_2\\  -(K-\frac{1}{m}G^2 + \lambda \Gamma)
&  -\frac{G}{m}\end{array}\right)
\AND
{F}_i=\left(\begin{array}{c}  0 \\ f_i \end{array}\right).
\]
Let $Z=(Z_1,\ldots,Z_N)^T$ and $ F=( F_1,\ldots,  F_N)^T$ represent now the state of the entire ring system, so that the ring
dynamics can be described in the following vector form
\begin{equation}\label{eq:Zsystem}
\dot{Z} = M Z-{F}(Z).
\end{equation}
Let $\tau=t$ and consider the system in extended phase space
\begin{equation}\label{eq:Zsystem-ext}
\dfrac{dZ}{dt} = MZ-F(Z) + H_{A_d}(\tau):=G(Z,\tau,A_d), \quad  \dfrac{d\tau}{dt} =1.
\end{equation}
where
\[
H_{A_d}(\tau) =\underbrace{(\underbrace{0,f_e(\tau),0,0},\ldots,\underbrace{0,f_e(\tau),0,0}).}_{\mbox{$N$ times}}
\]
Consider an equilibrium solution $Z_0$ of the unforced system~\eqref{eq:Zsystem}. One can show using the implicit function theorem that for small $2\pi/\omega_d$-periodic forcing,  if a non-resonance condition on the eigenvalues of the Jacobian at $Z_0$ is satisfied, there exists a $2\pi/\omega_d$-periodic solution $P_0(t)$ of the forced system~\eqref{eq:Zsystem-ext} passing near $Z_0$ in the extended phase space. See Chicone~\cite{Chicone-book} for the proof.
Suppose $Z_0$ has isotropy subgroup $\Sigma$. By unicity of the existence of $P_0(t)$ (from the implicit function theorem) and the uniformity of $H_{A_d}(\tau)$ along the $N$ gyroscopes, we can restrict~\eqref{eq:Zsystem} to $\Fix(\Sigma)$  and use the same implicit function theorem argument. Therefore, $P_0(t)$ lies in $\Fix(\Sigma)$.

Moreover, the monodromy matrix $M(2\pi/\omega_d)$ at $P_0(t)$ is obtained by solving the variational system
of~\eqref{eq:Zsystem-ext} at $P_0(t)$
\[
\dfrac{d\zeta}{dt} = dG(P_0(t))\zeta, \qquad \dfrac{d\xi}{dt}=0
\]
with $(\zeta,\xi)\in \R^{4N}\times\R$ and by the Floquet theorem, the spectrum of $M(2\pi/\omega_d)$ is given by the spectrum of $e^{(2\pi/\omega_d) dG(Z_0,0,0)}$ and the simple eigenvalue $+1$. These results are summarized in the next statement.
\begin{proposition}\label{prop:eq-per}
For sufficiently small forcing amplitude $A_d$, equilibrium solutions of~\eqref{eq:Zsystem} with isotropy subgroup $\Sigma$ are in one-to-one correspondence with $2\pi/\omega_d$-periodic solutions of~\eqref{eq:Zsystem-ext} with spatial symmetry group $\Sigma$. The Floquet exponents of the $2\pi/\omega_d$-periodic solution are given by the eigenvalues of the linearization of the corresponding equilibrium solution.
\end{proposition}

In the following sections, we look at the effect of coupling on the Hamiltonian structure. We show that in the unidirectional case, system~\eqref{eq:gen-coupling} is not Hamiltonian, while in the bidirectional case the system possesses a Hamiltonian structure.

\subsection{Unidirectional Coupling}
With  $N$ identical gyroscopes coupled unidirectionally,~\eqref{eq:gen-coupling} becomes
\begin{equation}\label{eq:unidir}
\left(\begin{array}{c}
\dot{q_i} \\ \dot{p}_i \end{array}\right) =
\left(\begin{array}{cc} -\frac{G}{m} & \frac{1}{m}I_2\\  -(K-\frac{1}{m}G^2 + \lambda \Gamma) &  -\frac{G}{m}\end{array}\right)\left(\begin{array}{c} q_i \\ p_i \end{array}\right)  + \left(\begin{array}{c} 0 \\ -f_i+ \lambda \Gamma q_{i+1} \end{array}\right).
\end{equation}
In a laboratory experiment, this type of coupling configuration could be realized by a microcircuit where the oscillations of any of the driving axes are processed electronically and input into the driving axis of the next unit in a cyclic manner. Although experimental works are important, we focus on the theoretical aspects of the rings dynamics in this paper. Thus, we now show proof that a ring of unidirectionally coupled gyroscopes does not possess Hamiltonian structure.

\begin{proposition}\label{prop:unidir}
The unidirectionally coupled gyroscopic system formulated using system~\Ref{eq:unidir} is not Hamiltonian with respect to the symplectic
structure given by
\begin{equation}\label{eq:Jmat}
J=\diag\underbrace{(J_4,\ldots,J_4)}_{N\;\mbox{times}},
\end{equation}
with $J_4 =  \left(\begin{array}{cc} 0 & I_2 \\-I_2 & 0 \end{array}\right)$.
\end{proposition}

\proof Consider~\eqref{eq:Zsystem} where
\begin{equation} \label{eq:linearSystem}
M=
\left(\begin{array}{cccccc}
M_1 & M_2 & 0 &\ldots&0\\
0 & M_1 & M_2&\ldots&\vdots \\
\vdots&0&\ddots&\ddots&\vdots\\
0&0&\ldots&\ddots&M_2\\
M_2 & 0 &\hdots&0& M_1
\end{array}\right)
\qquad \mbox{with}\qquad
M_2 = \left(\begin{array}{cc} 0 & 0 \\ \lambda\Gamma & 0  \end{array}\right).
\end{equation}

We can check directly that $M_1$ and $M_2$ are Hamiltonian matrices with respect to $J_4$.
That is, $M_j^{T} J_4 + J_4 M_j =0$ for $j=1,2$. By definition, the quadratic part of the
Hamiltonian function of the system can also be written in terms of some symmetric
matrix $S$ as
\[
H(Z) = \frac{1}{2}Z^T S Z,
\]
where $S=J^{-1}M$.

Therefore, if the unidirectionally coupled system in~\eqref{eq:unidir} were to admit a Hamiltonian
structure we should be able to find a matrix $S$ such that $S = -JM$ and $S=S^T$. After some computations, we get
\[
S=
\left(\begin{array}{cccccc}
S_1 & S_2 & 0 &\ldots&0\\
0 & S_1 & S_2&\ldots&\vdots \\
\vdots&0&\ddots&\ddots&\vdots\\
0&0&\ldots&\ddots&S_2\\
S_2 & 0 &\hdots&0& S_1
\end{array}\right).
\]
Since $S \neq S^T$, $M$ is not a Hamiltonian matrix and the coupled gyroscopic system formulated by~\eqref{eq:unidir}
is not a Hamiltonian system. \qed

We do not continue studying this case. Instead, we focus on the bidirectional case, which retains the Hamiltonian structure as it is shown in the next section.

\subsection{Bidirectional Coupling}
\label{sec:bidrectionalyCoupledGyroscopes}
For the bidirectional case, we use the appropriate coupling function so that the system~\eqref{eq:gen-coupling}
can be re-written as
\begin{equation}\label{eq:bidir}
\left(\begin{array}{c} \dot{ q_i} \\ \dot{ p}_i \end{array}\right) =
\left(\begin{array}{cc}
-\frac{G}{m} & \frac{1}{m}I_2\\  -(K-\frac{1}{m}G^2 + 2\lambda \Gamma) &  -\frac{G}{m}\end{array}\right)\left(\begin{array}{c} q_i \\
p_i \end{array}\right)  + \left(\begin{array}{c} 0 \\  \lambda \Gamma ( q_{i+1}+ q_{i-1}) -f_i \end{array}\right).
\end{equation}

This type of coupling can also be realized, in principle, electronically through a microcircuit as it was
described in the unidirectionally case. Additionally, bidirectional coupling could be easier to implement
in hardware by connecting the proof mass of adjacent gyroscopes through springs. Again, experimental
works are beyond the scope of the present manuscript. Instead, we show next that a bidirectionally
coupled gyroscope system possesses a Hamiltonian structure.

\begin{proposition}\label{prop:bidir}
The bidirectionally coupled gyroscopic system formulated through~\Ref{eq:bidir} is Hamiltonian with
respect to the symplectic structure given by
\[
J=\diag\underbrace{(J_4,\ldots,J_4)}_{N\;\mbox{times}},
\]
with $J_4 =  \left(\begin{array}{cc} 0 & I_2 \\-I_2 & 0 \end{array}\right)$.
\end{proposition}

\proof Consider the $N\times N$ cyclic permutation matrix
\begin{equation*}
C=\left(\begin{array}{cccccc}
0&1&0&0&\ldots&0\\
0&0&1&0&\ldots&0\\
\vdots&\vdots&\ddots&\ddots&\ldots&\vdots\\
\vdots&\ldots&&\ddots&\ddots&0\\
0&\ldots&\ldots&&0&1\\
1&0&\ldots&&\ldots&0
\end{array}\right).
\end{equation*}
Recall that the Kronecker product $\otimes$ of two matrices $A=[a_{ij}]$  of size $m\times n$ and $B$ of size
$p\times q$ is a $mp\times nq$ matrix defined by
\[
A\otimes B = [a_{ij} B].
\]
Let
\[
M = I_{N} \otimes M_1 + (C+C^{T}) \otimes M_2
\]
where $I_{N}$ is the $N$-dimensional identity matrix and $M_2$ is the same matrix as in the
unidirectional coupling case, see~\eqref{eq:linearSystem}. The matrix $M_1$ is, however,
slightly different
\[
M_1 = \left(\begin{array}{cc} -\frac{G}{m} & \frac{1}{m}I_2\\  -(K-\frac{1}{m}G^2 + 2\lambda \Gamma)
&  -\frac{G}{m}\end{array}\right).
\]

Once again direct calculations show that $M_1$ and $M_2$ are Hamiltonian matrices with respect to $J_4$, so that $M_j^{T} J_4 + J_4 M_j =0$ for $j=1,2$. We use again $Z_i = (q_i,p_i )^T$ to represent the position and momentum coordinates of the $i^{th}$ gyroscope and $Z=(Z_1,\ldots,Z_N)^{T}$ to represent the state of the entire ring at any time $t$. The governing equations for the bidirectionally coupled ring~\Ref{eq:bidir} can now be rewritten as
\begin{equation}\label{eq:hamiltonianSystem}
\dot Z = MZ- F(Z),
\end{equation}
where
\[
{F}(Z) = \left(\begin{array}{c} F_1\\ \vdots \\ F_N  \end{array}\right).
\]
Since $M_1$ and $M_2$ are Hamiltonian with respect to $J_4$, a direct calculation shows
that $M$ satisfies the condition $M^{T}J+JM=0$, and thus $M$ is a Hamiltonian matrix.
Finally, $F(Z)=J\nabla H_2(q,p)$, where $\displaystyle H_2(q,p)= \frac{1}{4}\sum_{i=1}^{N}\mu\left(q_{i1}^4+q_{i2}^4\right)$
and this completes the proof. \qed

We now complete the computation of the Hamiltonian function associated with the system in~\eqref{eq:hamiltonianSystem}.
Let $S = J^{-1}M$ so that the Hamiltonian function corresponding to the linear part of the
system is $H_0 = \frac{1}{2} Z^TSZ$. By the definition, $S$ has the form
\[
S = I_N \otimes S_1 + (C+C^{T})\otimes S_2,
\]
where $S_1 = J_4^T M_1$ and $S_2 = J_4^T M_2$. Note that $S=S^T$ and the corresponding linear Hamiltonian function is
\begin{equation*}
H(Z) = \frac{1}{2}Z^TSZ = \frac{1}{2} \sum_{i=1}^{N} Z_i^T S_1Z_i + \frac{1}{2} \sum_{i=1}^{N} (Z_{i+1}^T + Z_{i-1}^T)S_2Z_i.
\end{equation*}
The Hamiltonian of the complete $\D_N$-symmetric bidirectionally coupled system can now be expressed in terms of the position and momentum coordinates as
\begin{align*}
 H(q,p) ={}& \frac{1}{2}\sum_{i=1}^{N} -p_i^T\left(K-\frac{G^2}{m}+2\lambda\Gamma\right)q_i- q_i^T\frac{G}{m}q_i+p_i^T\frac{G}{m}p_i+q_i^T \frac{I_2}{m} p_i \\
 {}& -(q_{i+1}+q_{i-1})^T\lambda\Gamma q_i + H_2(q,p).
\end{align*}

\section{Linear Analysis at the Origin}\label{sec:linear-origin}
We begin the study of the linearized system near the origin starting with the $\D_N$ isotypic decomposition
of the tangent space. This leads to a block diagonal decomposition from which the eigenvalues are obtained
explicitly and their distribution is studied for all $N\in \N$. In particular, we determine for general $N$, a threshold condition for the origin to  lose spectral stability as the coupling parameter $\lambda$ is varied. The eigenvalue structure at the origin also enables us to determine the Lyapunov families of periodic orbits via the Equivariant Weinstein-Moser theorem.

\subsection{Isotypic Decomposition}\label{sec:isotypicDecomposition}
After transforming~\Ref{eq:coupled_gyro} into its Hamiltonian form, an additional simplification
is carried out by decomposing the system into its isotypic components~\cite{golubitsky1988groups}.
If we let
\[
{\cal K} = \left(\begin{array}{ccccccc}
0 & 1 & 0 & \ldots & \ldots & \ldots & 0 \\
1 & 0 & 0 & \ldots & \ldots & \ldots & 0 \\
0 & 0 & 0 & \ldots & \ldots & 0 & 1\\
\vdots&\vdots&0&\ldots &0&1 & 0\\
\vdots&\vdots&\vdots&\iddots&\iddots&\iddots&\vdots\\
\vdots&\vdots&0&\iddots&\iddots&&\vdots\\
0&0&1&0& \ldots &\ldots& 0\\
\end{array}\right).
\]
then we can write the generators of $\D_N = \langle \gamma,\kappa \rangle$ in $\mathbb{R}^{4N}$ in terms of the matrices $C$ and $\cal K$ as
\begin{equation}\label{eq:generators}
\gamma = C \otimes I_4 \AND
\kappa= {\cal K} \otimes I_4.
\end{equation}
The isotypic decomposition of $\C^N$ by $\langle C,{\cal K}\rangle$ is well-known, see~\cite{golubitsky1988groups}, and is given by
\[
\C^N = V_0 \oplus V_1 \oplus \cdots \oplus V_{N-1},
\]
where
\[
V_{j} = \C\{v_j\} \quad \mbox{with} \quad v_{j}=(v,\zeta^j v,\zeta^{2j} v,\ldots,\zeta^{(N-1)j}v)^T
\quad \textnormal{and} \quad \zeta = \exp\left({2\pi i/N}\right),
\]
for some $v\in \R$. Therefore, the isotypic decomposition of the complexified phase space is
\begin{equation}\label{eq:decCN4}
\left(\C^{N}\right)^{4} = V_0^{4} \oplus V_1^{4} \oplus \cdots \oplus V_{N-1}^{4}
\end{equation}
and
\[
{\cal M}_{j} := M\mid_{V_{j}} = M_1+(\zeta^{j}+\zeta^{(N-1)j})M_2 = M_1+(\zeta^{j}+\ov{\zeta}^{j})M_2 = M_1+2\cos(2\pi j/N) M_2,
\]
for $j=0,1,\ldots,(N-1)$. Note that if $j\neq 0, N/2$ and $N$ is even, ${\cal M}_{j} = {\cal M}_{(N-1)j}$.

We now verify that the basis of the decomposition formulated in~\eqref{eq:decCN4} is symplectic.
Let
\[
e_{1} = \left(\begin{array}{c} 1 \\ 0 \\ 0 \\ 0 \end{array}\right), \quad
e_{2} = \left(\begin{array}{c} 0 \\ 1 \\ 0 \\ 0 \end{array}\right), \quad
e_{3} = \left(\begin{array}{c} 0 \\ 0 \\ 1 \\ 0 \end{array}\right), \quad
e_{4} = \left(\begin{array}{c} 0 \\ 0 \\ 0 \\ 1 \end{array}\right),
\]
and define
\[
v_{ji} = (e_i,\zeta^j e_i,\zeta^{2j} e_i,\ldots,\zeta^{(N-1)j}e_i)^T,
\]
for $i=1,2,3,4$, and $j=0,\ldots,N-1$. We need to verify that for the symplectic form
$\omega(u,v)=u^T Jv$ with $u,v\in \C^{4N}$ and $J$ given by~(\ref{eq:Jmat}),
$\omega(v_{ji},v_{k\ell})=0$ holds for any pair $v_{ji},v_{k\ell}$ in the basis of $\C^{4N}$. We have
\[
\omega(v_{ji},v_{k\ell}) = v_{ji}^{T} J v_{\ell k} = \sum_{m=0}^{N-1} \zeta^{m(j+\ell)} e_i^{T} J_4 e_k=e_i^{T} J_4 e_k  \sum_{m=0}^{N-1} \zeta^{m(j+\ell)},
\]
and note that
\[
\sum_{m=0}^{N-1} \zeta^{m(j+\ell)} = 0,
\]
for any combination of $j,\ell$. The corresponding real symplectic transition matrix $P$ is constructed using the normalized real and imaginary parts of the vectors $v_{ji}$ for complex vectors and just the normalized $v_{ji}$ if it is real. For complex vector $v_{ji}$, let $\Im_{ji}$ and $\Re_{ji}$ denote their imaginary and real parts, respectively. Furthermore, we denote a normalized vectors by $\tilde{\cdot}$. For $N$ odd, the real symplectic transition matrix is
\[
P=\left[\tilde{v}_{01},\ldots,\tilde{v}_{04},\widetilde{\Im}_{11},\ldots,\widetilde{\Im}_{14},\widetilde{\Re}_{11},\ldots,\widetilde{\Re}_{14},\ldots,\widetilde{\Im}_{\left\lfloor {N}/{2}\right\rfloor1},\ldots,\widetilde{\Im}_{\left\lfloor {N}/{2}\right\rfloor4},\widetilde{\Re}_{\left\lfloor {N}/{2}\right\rfloor1},\ldots,\widetilde{\Re}_{\left\lfloor {N}/{2}\right\rfloor4}\right].
\]
Similarly, the corresponding real symplectic matrix for $N$ even is
\begin{align*}
P={}&\left[\tilde{v}_{01},\ldots,\tilde{v}_{04},\widetilde{\Im}_{11},\ldots,\widetilde{\Im}_{14},\widetilde{\Re}_{11},\ldots,\widetilde{\Re}_{14},\ldots, \widetilde{\Im}_{\left({N}/{2}-1\right)1},\ldots,\widetilde{\Im}_{\left({N}/{2}-1\right)4},\widetilde{\Re}_{\left({N}/{2}-1\right)1},\ldots,\right.\\
&\quad\widetilde{\Re}_{\left({N}/{2}-1\right)4},\tilde{v}_{\left({N}/{2}\right)1},\ldots,\tilde{v}_{\left({N}/{2}\right)4}\Big].
\end{align*}

Applying $P$ to the linear part of~\eqref{eq:hamiltonianSystem} we obtain
\[
P^{-1}MP = \mathcal M,
\]
where for $N$ odd,
\begin{equation*}
\mathcal M = \diag\left({\cal M}_0,{\cal M}_{1},{\cal M}_1,\ldots,{\cal M}_{\left\lfloor {N}/{2}\right\rfloor},{\cal M}_{\left\lfloor {N}/{2}\right\rfloor}\right)
\end{equation*}
and for $N$ even
\begin{equation*}
\mathcal M = \diag\left({\cal M}_0,{\cal M}_{1},{\cal M}_1,\ldots,{\cal M}_{{N}/{2}-1},{\cal M}_{{N}/{2}-1},{\cal M}_{{N}/{2}}\right)
\end{equation*}
Because $\mathcal M$ is Hamiltonian, every  ${\cal M}_j$ block is also Hamiltonian with respect to $J_4$. Thus, there is a symmetric matrix $S_\mathcal M$ such that $JS_\mathcal M = \mathcal M$. When $N$ is odd, the corresponding symmetric matrix is
\[
S_{\mathcal M} = \diag\left(J_4^{T} {\cal M}_0,J_4^{T} {\cal M}_1,J_4^{T} {\cal M}_1,\ldots,J_4^{T} {\cal M}_{\left\lfloor {N}/{2}\right\rfloor}, J_4^{T} {\cal M}_{\left\lfloor {N}/{2}\right\rfloor}\right).
\]
Similarly, when $N$ is even, the symmetric matrix is
\[
S_{\mathcal M} = \diag\left(J_4^{T} {\cal M}_0,J_4^{T} {\cal M}_1,J_4^{T} {\cal M}_1,\ldots,J_4^{T}{\cal M}_{{N}/{2}-1},J_4^{T}{\cal M}_{{N}/{2}-1},J_4^{T} {\cal M}_{{N}/{2}}\right).
\]
Let $S_\mathcal M^j=J_4^T {\cal M}_j$, where $j=0,\ldots,\left\lfloor {N}/{2}\right\rfloor$. Using the change of variables $Y=PU$, where $U=(U_1,\ldots,U_N)$, we obtain the quadratic Hamiltonian function
\[
H_0^{odd}(U) =\frac{1}{2}U^T S_\mathcal M U = \dfrac{1}{2} U_1^{T} S_{\mathcal M}^{0} U_1 + \dfrac{1}{2}\sum_{m=1}^{\left\lfloor {N}/{2}\right\rfloor} U_{m+1}^{T} S_{\mathcal M}^{m} U_{m+1} + U_{m+2}^{T} S_{\mathcal M}^{m} U_{m+2}
\]
for $N$ odd and for $N$ even, we have
\begin{align*}
H_0^{even}(U) ={}&\frac{1}{2}U^T S_{\mathcal M} U = \dfrac{1}{2} U_1^{T} S_{\mathcal M}^{0} U_1+ \dfrac{1}{2} U_{{N}/{2}+1}^{T} S_{\mathcal M}^{N/2} U_{{N}/{2}+1}+\dfrac{1}{2}  \sum_{m=1}^{{N}/{2}-1} U_{m+1}^{T} S_{\mathcal M}^{m} U_{m+1} \\
{}&+ U_{m+2}^{T} S_{\mathcal M}^{m} U_{m+2}.
\end{align*}

\subsection{Eigenvalues}
\label{sec:eigenvalues}
We now calculate the eigenvalues of the matrix $\mathcal M$, which are required to further simplify
the Hamiltonian system into normal form. Since the linear system is in block diagonal form, the
eigenvalues are the same as the combined eigenvalues of all the $\mathcal M_j$ blocks. In general,
each block can be written as
\[
{\cal M}_{j} = M_1 + 2\cos{(2\pi j/N)} M_2 =
\left(\begin{array}{cc} -\frac{G}{m} & \frac{1}{m}I_2\\  -(K-\frac{1}{m}G^2 + 2\lambda \Gamma(1-\cos{(2\pi j/N)})
&  -\frac{G}{m}\end{array}\right).
\]
Of the four eigenvalues, two of them have the form
\[
     \rho_{j}^{\pm} = {1 \over \sqrt{m}}\sqrt{-\left(\kappa + 2m\Omega^2 + \lambda \left(1-\cos{2 \pi j \over N} \right)\right)\pm \sqrt{s_j}},
\]
where $s_j=4m\Omega^2(\kappa+m\Omega^2+\lambda(1-\cos{(2\pi j/N)}))+\lambda^2 (1-\cos{(2\pi j/N)})^2$.
The other two eigenvalues are $-\rho_{j}^{\pm}$. It is straightforward to check that because $\kappa>0$, the eigenvalue $\rho_{j}^{-}$ is purely imaginary for all $\lambda\in\R$. Observe that $\rho_{j}^{+}=0$ if and only if $\kappa+2\lambda(1-\cos(2\pi j/N))=0$, that is,
\begin{equation}\label{param:crit}
\lambda_{j}^{*} = \frac{-\kappa}{2(1-\cos(2\pi j/N))}.
\end{equation}
This result implies that $\lambda_j^*$ is maximum when $j=\left\lfloor{N}/{2}\right\rfloor$. For $N$ even,
$\lambda_{\left\lfloor{N}/{2}\right\rfloor}^{*}=-\kappa/4$ and for $N$ odd $\lambda_{\left\lfloor{N}/{2}\right\rfloor}^{*}$ takes its smallest value for $N=3$ at $-\kappa/3$ and converges to $-\kappa/4$ as $N\to \infty$.

One can easily check that as $\lambda$ increases through $\lambda_{j}^*$, $\rho_{j}^+$ changes from real
to purely imaginary. Thus, for $\lambda>\lambda_{\left\lfloor{N}/{2}\right\rfloor}^*$ all eigenvalues
are purely imaginary. Recalling that an equilibrium is {\em spectrally stable} if all the eigenvalues of
the linearization of the equilibrium are on the imaginary axis, we then arrive at the following threshold
condition for stability.

\begin{proposition}\label{prop:crit-eig}
For $\lambda>\lambda_{\left\lfloor{N}/2\right\rfloor}^*$, the equilibrium at the origin is spectrally stable and unstable for $\lambda<\lambda_{\left\lfloor{N}/2\right\rfloor}^*$. Moreover, if $\lambda<-\kappa/2$ then $\rho_{j}^{+}\in \R$ for all $j=0,\ldots,\left\lfloor{N}/{2}\right\rfloor$.
\end{proposition}

For $N$ even, $\lambda_{N/2}^*$ is the threshold value for a bifurcation from the ${\cal M}_{N/2}$ block
and so a single pair of eigenvalues crosses the origin. This leads to a $\Z_2$ symmetry-breaking bifurcation.
That is, a pitchfork bifurcation. For $N$ odd, $\lambda_{\left\lfloor{N}/{2}\right\rfloor}^*$ is the threshold value for a bifurcation from the two ${\cal M}_{\left\lfloor{N}/{2}\right\rfloor}$ blocks and thus a double pair of eigenvalues crosses the origin. Therefore, this is a $\D_N$ symmetry-breaking bifurcation where a group orbit of equilibria with isotropy subgroup $\Z_2$ (and its conjugates) bifurcate from the origin via the Equivariant Branching Lemma. A more complete picture of the nature of the bifurcations is obtained later on via normal form analysis and the equivariant splitting lemma.

Observe that the term $f_e(t)$ from~\Ref{eq:coupled_gyro} does not play a role in the linear part of the system and so it does not affect the calculation of the eigenvalues. Consequently, excluding the damping terms for right now, $\lambda_{N/2}^*$ represents an analytical approximation to the critical coupling
strength that leads the gyroscopes into complete synchronization when the periodic forcing term $f_e(t)$
is added~\cite{vu2011drive}.

We now look at the distribution of eigenvalues on the imaginary axis.
\begin{proposition}\label{prop:dist-eig}
Let $j$ increase from $0$ to $\left\lfloor{N}/{2}\right\rfloor$. Then, for $\lambda<0$, $\im(\rho_{j}^{\pm})$ decreases as a function of $j$ and for $\lambda>0$, $\im(\rho_{j}^{\pm})$ increases as a function of $j$.
\end{proposition}

\proof Treat $j$ as a continuous variable and take the derivative of $\im(\rho_{j}^{-})$ to obtain
\begin{equation}\label{eq:deriv-imag}
\frac{\pi \lambda}{N\im(\rho_{j}^{-})}\sin(2\pi j/N)\left(1+\frac{1}{2\sqrt{s_j}}(4m\Omega^2+2\lambda(1-\cos(2\pi j/N)))\right).
\end{equation}
For $\lambda>0$, the derivative~(\ref{eq:deriv-imag}) is positive for $j\in [0,\left\lfloor{N}/{2}\right\rfloor]$. For
$\lambda<0$, ~(\ref{eq:deriv-imag}) is negative because $|\lambda(1-\cos(2\pi j/N))/\sqrt{s_j}|<1$.
For $\rho_{j}^{+}$, the derivative of $\im(\rho_{j}^{+})$ is
\[
\frac{\pi \lambda}{N\im(\rho_{j}^{-})}\sin(2\pi j/N)\left(1-\frac{1}{2\sqrt{s_j}}(4m\Omega^2+2\lambda(1-\cos(2\pi j/N)))\right).
\]
where the domain of $j$ is shrinked correspondingly if $\lambda<\lambda_{\left\lfloor{N}/{2}\right\rfloor}^{*}$. Inspection of $s_j$ shows that
\[
\frac{1}{2\sqrt{s_j}}(4m\Omega^2+2\lambda(1-\cos(2\pi j/N)))<1.
\]
and so the derivative has the sign of $\lambda$ as in the $\rho_{j}^{-}$ case. \qed

From Proposition~\ref{prop:dist-eig}, the purely imaginary eigenvalues $\rho_{j}^{\pm}$ are distributed monotonically and do not intersect for all $\lambda$.

Note that for $\lambda<0$, $\im(\rho_{j}^{-})>\im(\rho_{j}^{+})$ for all $j=0,\ldots,\left\lfloor{N}/{2}\right\rfloor$ and from  Proposition~\ref{prop:dist-eig}, $\im(\rho_{0}^{-})>\im(\rho_{j}^{-})$ for all $j=1,\ldots,\left\lfloor{N}/{2}\right\rfloor$. Because $\rho_{0}^{-}$ is independent of $N$, $\rho_{0}^{-}$ is an upper bound for all purely imaginary eigenvalues, for all $N\in \N$.  For $\lambda>0$, Proposition~\ref{prop:dist-eig} shows $\im(\rho_{\left\lfloor{N}/{2}\right\rfloor}^{-})>\im(\rho_{j}^{-})$ for all $j$ and note that the value $\im(\rho_{\left\lfloor{N}/{2}\right\rfloor}^{-})$ is bounded above by the following constant:
\[
{1 \over \sqrt{m}}\sqrt{\left(\kappa + 2m\Omega^2\right)+\sqrt{4m\Omega^{2}(\kappa+m\Omega^2)}}.
\]
Because of the upper bounds on the purely imaginary eigenvalues, the distance between nearby $\rho_{j}^{-}$ shrinks as $N$ increases. The same is true for $\rho_{j}^{+}$ as long as some of them are purely imaginary.

\subsection{Lyapunov Families} \label{sec:stablity}
In this section, we further the study of the local dynamics in the neighborhood of the equilibrium at the origin.
We show the existence of families of symmetric periodic orbits near the origin using the Equivariant
Weinstein-Moser (EWM) theorem, see Montaldi {\em et al.}~\cite{Montaldietal1988}.
To apply the Equivariant Weinstein-Moser theorem two conditions must be satisfied:
\begin{description}
   \item[(H1)] $D^2 H(p)$ must be a nondegenerate quadratic form,
   \item[(H2)] $D^2 H(p)|_{V_{\nu}}$ is positive definite,
\end{description}
where $V_{\nu}$ is the resonance subspace of the eigenvalue $\nu$ of the linearization at the origin.
Condition (H1) is satisfied at all values of $\lambda$ for which there are no zero eigenvalues. Condition
(H2) is satisfied for all purely imaginary eigenvalues $\nu$.

\begin{theorem}
For each eigenvalue $\rho_{j}^{\pm}\in i\R$ of the diagonal block ${\cal M}_{j}$, there exists at least one near $2\pi/|\rho_{j}^{\pm}|$-periodic solution for each energy level close to $H(p)$ with spatio-temporal isotropy subgroup
\begin{enumerate}
  \item $\Z_2(\kappa)$ ($N$-odd) or $\Z_2(\kappa)\times\Z_2^{c}$ ($N$-even),
  \item $\Z_2(\kappa,\pi)$ ($N$-odd) or $\Z_2(\kappa,\pi)\times\Z_2^{c}$,  $\Z_2(\kappa,\pi)\times\Z_2^{c}$ ($N$-even),
  \item $\Z_{N/gcd(N,j)}(\gamma^{j},2\pi j/N)$
\end{enumerate}
where $\Z_2^{c}=\Z_2(\gamma^{N/2},\pi)$. The first two are standing waves and the last one is a
discrete rotating wave.
\end{theorem}

\proof Suppose that ${\cal M}_{j}$ has only purely imaginary eigenvalues $\pm \rho_{j}^{\pm}$. Then,
\[
V_{j}^{4} = V_{j,\mu_{j}^{+}} \oplus V_{j,\mu_{j}^{-}},
\]
and condition (H2) is satisfied for both subspaces. From the Equivariant Weinstein-Moser theorem (and following remarks) and as shown in section 7 of~\cite{Montaldietal1988}, for each energy level near $H(p)$, at least one periodic solution of period $2\pi/|\rho^{\pm}|$ exists with symmetry corresponding exactly to one of the three (conjugacy classes of) isotropy subgroups of $\D_n\times \Sone$. One subgroup, $\tilde{\Z}_N$, is cyclic of order $N$ and it represents a rotating wave in which all gyroscopes oscillate with the same wave form and same amplitude but with phase shifts of $2\pi/N$ from one to the next. The other subgroups are isomorphic to $\Z_2$ (or $\Z_2 \oplus \Z_2$ when $N$ is even), but with subtle differences depending on wether $N=2$ (mod 4) or $N=0$ (mod 4). In either case, these two subgroups represent standing waves.
 \qed

In fact, because each of these three conjugacy classes of isotropy subgroups in $\D_n\times \Sone$ has fixed point subspace of dimension two, there exists a $C^{\infty}$ two-dimensional manifold passing through $p$ foliated by periodic solutions with periods near $2\pi/|\rho^{\pm}|$ and corresponding symmetry groups described above. The tangent space of the submanifold is tangent to the fixed-point subspaces $\Fix(\Sigma)$.
Because the homogeneous equilibrium $p$ is unstable for $\lambda<\lambda_{\left\lfloor{N}/{2}\right\rfloor}^{*}$, the periodic solutions near $p$ can only be stable for $\lambda>\lambda_{\left\lfloor{N}/{2}\right\rfloor}^{*}$.

\section{Normal Form Analysis}

\subsection{Linear Normal Form}
\label{sec:linearAnalysis}
In section~\ref{sec:bidrectionalyCoupledGyroscopes}, we show that a system of $N$ bidirectionally coupled gyroscopes is Hamiltonian. The next step in the analysis is to determine the normal form of the coupled system at the bifurcation point given by~\eqref{param:crit}. We begin the calculations by finding a symplectic matrix $Q$ to transform the linear part of the system into normal form~\cite{meyer2008introduction}. Given the block diagonal structure of $\mathcal M$, we can construct symplectic transition matrices $Q_j$ corresponding to each $\mathcal M_j$ and combine them to form $Q$ as desired. As shown in section~\ref{sec:eigenvalues}, when $j\neq\left\lfloor{N}/{2}\right\rfloor$, there are two pairs of purely imaginary eigenvalues for each $\mathcal M_j$. However, when $j=\left\lfloor{N}/{2}\right\rfloor$, the eigenvalues of $Q_j$ consist of a pair of zeroes and a pair of purely imaginary eigenvalues. Thus, we consider the two cases separately.

For $j\neq\left\lfloor{N}/{2}\right\rfloor$, we can apply the method outlined in~\cite{burgoyne1974normal} to obtain the corresponding symplectic transformation $Q_j$. To begin, we write the eigenvalues in complex form as
\[
     i\nu_{j}^{\mp} =
     {\frac{i}{ \sqrt{m}}}\sqrt{\left(\kappa + 2m\Omega^2 +
                             \lambda \left(1-\cos{2 \pi j \over N} \right)\right)\mp \sqrt{s_j}},
\] where $i=\sqrt{-1}$. Setting $\lambda=\lambda^*_{\left\lfloor{N}/{2}\right\rfloor}$,  $\displaystyle \omega=2\lambda^*_{\left\lfloor{N}/{2}\right\rfloor}\left(1-\cos{\frac{2\pi j}{N}}\right)$ and $q=\sqrt{16m\Omega^2\kappa+(\kappa\omega-4m\Omega^2)^2}$, the symplectic matrix $Q_j$ is
\begin{equation}\label{eq:qJnotEqualNover2}
\scriptsize
Q_j=
\left(
\begin{array}{cccc}
\sqrt{\frac{{\nu^-_j}}{c_1}}\left(\frac{4\Omega}{-\kappa\omega+q}\right)&0&0&\sqrt{\frac{{\nu^+_j}}{2c_{2}}}\left(\frac{-4\Omega}{\kappa\omega+q}\right)
\\
0&\frac{1}{\sqrt{2c_2\nu^+_j}}\left(\frac{4\Omega^2}{\kappa\omega+q}+\frac{1}{m}\right)&\frac{1}{\sqrt{c_1\nu^-_j}}\left(\frac{-4\Omega^2}{-\kappa\omega+q}+\frac{1}{m}\right)&0
\\
0&\frac{1}{\sqrt{2c_2\nu^+_j}}\left(\frac{4(\kappa+m\Omega^2-\kappa\omega)\Omega}{\kappa\omega+q}+\Omega\right)&\frac{1}{\sqrt{c_1\nu^-_j}}\left(\frac{4(-\kappa-m\Omega^2+\kappa\omega)\Omega}{-\kappa\omega+q}-\Omega\right)&0
\\
\sqrt{\frac{{\nu^-_j}}{c_1}}&0&0&\sqrt{\frac{{\nu^+_j}}{2c_{2}}}
\end{array}
\right),
\end{equation}
where $\displaystyle c_1=\frac{q^2-(4m\Omega^2+\kappa\omega)q}{(-\kappa\omega+q)^2m}$ and $\displaystyle c_2={\frac{q^2+(4m\Omega^2+\kappa\omega)q}{(\kappa\omega+q)^2m}}$. Applying $Q_j$ to $\mathcal{M}_j$, we get
\[
\mathsf M_j= Q_j^{-1}\mathcal M_jQ_j=
\left(\begin{array}{cccc}
0&0&\nu^-_j&0
\\
0&0&0&\nu^+_j
\\
-\nu^-_j&0&0&0
\\
0&-\nu^+_j&0&0
\end{array}\right).
\]

Consider the other case when $j=\left\lfloor{N}/{2}\right\rfloor$,  the eigenvalues of $\mathcal M_j$ consist of
\begin{equation*}
\pm i\psi=\pm i\sqrt{\frac{ (\kappa+4m\Omega^2)}{m}},
\end{equation*} and a pair of zero eigenvalues. Thus, the corresponding symplectic transition matrix is
\begin{equation}\label{eq:qJEqualsNover2}
Q_j=
\left(
\begin{array}{cccc}
0&0&-\frac{2\Omega m^{1\over 4}}{(\kappa+4m\Omega^2)^{\frac{3}{4}}}&\sqrt{\frac{{\kappa}}{m(\kappa+4m\Omega^2)}}
\\
\frac{1}{(m(\kappa+4m\Omega^2))^{\frac{1}{4}}}&\frac{2\Omega\sqrt{m}}{\sqrt{\kappa(\kappa+4m\Omega^2)}}&0&0
\\
\frac{m^{\frac{3}{4}}\Omega}{(\kappa+4m\Omega^2)^{\frac{1}{4}}}&-\frac{\sqrt{m}(\kappa+2m\Omega^2)}{\sqrt{\kappa(\kappa+4m\Omega^2)}}&0&0
\\
0&0&\frac{(\kappa+2m\Omega^2)m^{\frac{1}{4}}}{(\kappa+4m\Omega^2)^{\frac{3}{4}}}&\sqrt{\frac{m\kappa\Omega^2}{\kappa+4m\Omega^2}}
\end{array}
\right).
\end{equation}
Applying $Q_j$ to $\mathcal{M}_j$, we get
\[
\mathsf M_j= Q_j^{-1}\mathcal M_jQ_j=
\left(\begin{array}{cccc}
0&0&\psi&0\\
0&0&0&0\\
-\psi&0&0&0\\
0&-1&0&0
\end{array}\right).
\]
Thus, when $N$ is odd, we construct the overall symplectic transition matrix as
\begin{equation*}
Q=\operatorname {diag}(Q_0,Q_1,Q_1,\ldots,Q_{\left\lfloor{N}/{2}\right\rfloor},Q_{\left\lfloor{N}/{2}\right\rfloor})
\end{equation*} and the linear part of the gyroscopic system becomes
\begin{equation*}
\mathsf M=Q^{-1}\mathcal M Q=\operatorname{diag}(\mathsf M_0,\mathsf M_1,\mathsf M_1,\ldots,\mathsf M_{\left\lfloor{N}/{2}\right\rfloor},\mathsf M_{\left\lfloor{N}/{2}\right\rfloor}).
\end{equation*} Similarly, when $N$ is even, we construct the overall symplectic matrix as
\begin{equation*}
Q=\operatorname {diag}\left(Q_0,Q_1,Q_1,\ldots,Q_{{N}/{2}-1},Q_{{N}/{2}-1},Q_{\left\lfloor{N}/{2}\right\rfloor}\right).
\end{equation*}and the linear part of the gyroscopic system becomes
\begin{equation*}
\mathsf M=Q^{-1}\mathcal M Q=\operatorname{diag}\left(\mathsf M_0,\mathsf M_1,\mathsf M_1,\ldots,\mathsf M_{{N}/{2}-1},\mathsf M_{{N}/{2}-1},\mathsf M_{\left\lfloor{N}/{2}\right\rfloor}\right).
\end{equation*}

Let $$X=(X_0,X_1,Y_1,\ldots,X_{\left\lfloor{N}/{2}\right\rfloor-1},Y_{\left\lfloor{N}/{2}\right\rfloor-1},X_{\left\lfloor{N}/{2}\right\rfloor},Y_{\left\lfloor{N}/{2}\right\rfloor}),$$
where $X_i=(x_{i1},x_{i2},x_{i3},x_{i4})$ and $Y_i=(y_{i1},y_{i2},y_{i3},y_{i4})$ with $Y_{\left\lfloor{N}/{2}\right\rfloor}=0$ for $N$ even, represent the new state coordinates
for the entire ring. Under the transformation $U= QX$, the Hamiltonian function can now be written as
\begin{equation}\label{eq:hamiltonianNormalFormDegreeTwo}
H(X)=\tilde H_0(X)+H_2(X),
\end{equation}
where $\tilde H_0(X)$ and $H_2(X)$ represent  polynomials of degree two and four, respectively and $\tilde \dot$ denotes a function already in normal form. With $\tilde H_0(X)$ already in normal form, it can be written as
\begin{equation*}
\tilde H_0(X)= \frac{1}{2}X^TJ^{-1}\mathsf MX.
\end{equation*}
When $N$ is odd, the linear normal form is
\begin{align}\label{eq:degreeTwoHamiltonianNormalFormOdd}
\begin{split}
\tilde H_0^{odd}(X)={}&\frac{1}{2}\left(X_0^TJ^{-1}_4\mathsf M_0X_0+X_1^TJ^{-1}_4\mathsf M_1X_1+Y_1^TJ^{-1}_4\mathsf M_1Y_1+\cdots+\right.
\\
&\quad \left.X_{\left\lfloor{N}/{2}\right\rfloor}^TJ^{-1}_4\mathsf M_{\left\lfloor{N}/{2}\right\rfloor}X_{\left\lfloor{N}/{2}\right\rfloor}+Y_{\left\lfloor{N}/{2}\right\rfloor}^TJ^{-1}_4\mathsf M_{\left\lfloor{N}/{2}\right\rfloor}Y_{\left\lfloor{N}/{2}\right\rfloor}\right)\\
={}&\frac{1}{2}\Big(\nu_0^{-}\left(x_{01}^2+x_{03}^2\right)+\nu_0^{+}\left(x_{02}^2+x_{04}^2\right)+\nu_1^{-}\left(x_{11}^2+x_{13}^2\right)+\nu_1^{+}\left(x_{12}^2+x_{14}^2\right)\\
&+\nu_1^{-}\left(y_{11}^2+y_{13}^2\right)+\nu_1^{+}\left(y_{12}^2+y_{14}^2\right)+\cdots+\psi\left(x_{\left\lfloor{N}/{2}\right\rfloor 1}^2+x_{\left\lfloor{N}/{2}\right\rfloor 3}^2\right)
+x_{\left\lfloor{N}/{2}\right\rfloor 2}^2
\\
&\left.+\psi\left(y_{\left\lfloor{N}/{2}\right\rfloor 1}^2+y_{\left\lfloor{N}/{2}\right\rfloor 3}^2\right)+y_{\left\lfloor{N}/{2}\right\rfloor 2}^2\right).
\end{split}
\end{align}Likewise, when $N$ is even, the linear normal form is
\begin{align}\label{eq:degreeTwoHamiltonianNormalFormEven}
\begin{split}
\tilde H_0^{even}(X)={}&\frac{1}{2}\Big(X_0^TJ^{-1}_4\mathsf M_0X_0+X_1^TJ^{-1}_4\mathsf M_1X_1+Y_1^TJ^{-1}_4\mathsf M_1Y_1+\cdots +\\
&\quad X_{\left\lfloor{N}/{2}\right\rfloor -1}^TJ^{-1}_4\mathsf M_{{N}/{2}-2}X_{\left\lfloor{N}/{2}\right\rfloor-1}+Y_{\left\lfloor{N}/{2}\right\rfloor -1}^TJ^{-1}_4\mathsf M_{{N}/{2}-1}Y_{\left\lfloor{N}/{2}\right\rfloor -1}+\\
&\quad X_{\left\lfloor{N}/{2}\right\rfloor}^TJ^{-1}_4\mathsf M_{{N}/{2}}X_{\left\lfloor{N}/{2}\right\rfloor}\Big)\\
={}&\frac{1}{2}\Big(\nu_0^{-}\left(x_{01}^2+x_{03}^2\right)+\nu_0^{+}\left(x_{02}^2+x_{04}^2\right)+\nu_1^{-}\left(x_{11}^2+x_{13}^2\right)+\nu_1^{+}\left(x_{12}^2+x_{14}^2\right) +
\\
&\quad \nu_1^{-}\left(y_{11}^2+y_{13}^2\right)+\nu_1^{+}\left(y_{12}^2+y_{14}^2\right)+\cdots+\nu_{{N}/{2}-1}^{-}\left(x_{({N}/{2}-1) 1}^2+x_{({N}/{2}-1) 3}^2\right)+
\\
&\quad \nu_{{N}/{2}-1}^{+}\left(x_{({N}/{2}-1) 2}^2+x_{({N}/{2}-1) 4}^2\right)+\nu_{{N}/{2}-1}^{-}\left(y_{({N}/{2}-1) 1}^2+y_{({N}/{2}-1) 3}^2\right)+
\\
&\quad \nu_{{N}/{2}-1}^{+}\left(y_{({N}/{2}-1) 2}^2+y_{({N}/{2}-1) 4}^2\right)+x_{({N}/{2}) 2}^2+\psi\left(x_{({N}/{2}) 1}^2+x_{({N}/{2}) 3}^2\right)\Big).
\end{split}
\end{align}
The expression for $H_2$ is too long to be reproduced here, but it is simplified through normal form methods in section~\ref{sec:nonlinearNormalForm} and written explicitly for the $\D_3$ symmetric system as a case study in section~\ref{sec:bifurcationParameter}.

For $\lambda>\lambda_{\lfloor{N}/{2}\rfloor}^{*}$, the block ${\cal M}_{\lfloor{N}/{2}\rfloor}$ has purely imaginary eigenvalues and from the results of this section, the Hamiltonian of the linear normal form would contain instead the term
\[
\begin{array}{l}
\nu_{\lfloor{N}/{2}\rfloor}^{-}(x_{(\lfloor{N}/{2}\rfloor)1}^2+x_{(\lfloor{N}/{2}\rfloor)3}^{2})+\nu_{\lfloor{N}/{2}\rfloor}^{+}(x_{(\lfloor{N}/{2}\rfloor)2}^2+x_{(\lfloor{N}/{2}\rfloor)4}^{2})+\nu_{\lfloor{N}/{2}\rfloor}^{-}(y_{(\lfloor{N}/{2}\rfloor)1}^2
+y_{(\lfloor{N}/{2}\rfloor)3}^{2})\\+\nu_{\lfloor{N}/{2}\rfloor}^{+}(y_{(\lfloor{N}/{2}\rfloor)1}^2+y_{(\lfloor{N}/{2}\rfloor)3}^{2}).
\end{array}
\]
where $\nu_{\lfloor{N}/{2}\rfloor}^{\pm}>0$. We have the following result.

\begin{proposition}
For $\lambda>\lambda_{\lfloor{N}/{2}\rfloor}^{*}$, the equilibrium solution at the origin is locally Lyapunov stable.
\end{proposition}

\proof From the form of the quadratic Hamiltonian $\tilde{H}_{0}$ (for both $N$ odd and even), the origin is a strict local minimum of $H$ and so by Dirichlet's theorem it is locally Lyapunov stable. \qed


\subsection{Nonlinear Normal Form}
\label{sec:nonlinearNormalForm}
The normal form theory for $\Gamma$-symmetric Hamiltonian systems states that one can find formal symplectic changes of coordinates such that $H_j$, the transformed Hamiltonian function at every homogeneous degree is $\Gamma$-invariant, and $H_{j}$ satisfies
\begin{equation}\label{cond:nf}
H_{j}(e^{tL^{T}}X) = H_{j}(X),
\end{equation}
where $L$ is the linearization at the equilibrium point, see Montaldi {\em et al}~\cite{Montaldietal1990}. We denote by ${\bf S}$ the closure of the group generated by $e^{tL^{T}}$. This means $H_{j}$ is $\Gamma\times{\bf S}$-invariant at every order of transformation. In our case, $\Gamma=\D_N\times \Z_2$ and $L={\cal M}$. It is a straightforward calculation using the diagonal blocks of ${\cal M}$ in linear normal form from Section~\ref{sec:linearAnalysis} to show that
${\bf S}=\ov{e^{t{\cal M}}} \simeq \T^{m}$ where $m=2\left\lfloor{N}/{2}\right\rfloor+1$ for $N$ odd and $m=N$ for $N$ even and so $H_{j}$ commutes with $\D_N\times\Z_2\times \T^{m}$.

The symplectic matrix $P$ that transforms ${\cal M}$ into block diagonal form also transforms the generators in~\eqref{eq:generators} into the new coordinates. When $N$ is odd, the generators can be written as
\begin{equation*}
\tilde{\gamma}=P^{-1}\gamma P=\left(I_4, I_4\otimes \mathcal R\left(\omega_1\right),\ldots,I_4\otimes \mathcal R\left (\omega_{\left\lfloor{N}/{2}\right\rfloor}\right)\right)
\end{equation*}and
\begin{equation*}
\tilde{\kappa}=P^{-1}\kappa P=\left(I_4, I_4\otimes \mathcal S\left(\omega_1\right),\ldots,I_4\otimes \mathcal S\left (\omega_{\left\lfloor{N}/{2}\right\rfloor}\right)\right),
\end{equation*}
where $$\omega_j=\frac{2\pi j}{N},\quad\mathcal R(\theta)=\left(\begin{array}{cc}\cos\theta&-\sin\theta\\\sin\theta&\cos\theta\end{array}\right),\AND\mathcal S(\theta)=\left(\begin{array}{cc}-\cos\theta&\sin\theta\\\sin\theta&\cos\theta\end{array}\right).$$
Since the calculations for the invariants are simpler to perform in complex coordinates, we identify $X=(X_0,X_1,Y_1,\ldots,X_{\left\lfloor{N}/{2}\right\rfloor},Y_{\left\lfloor{N}/{2}\right\rfloor})$, (defined as before) in $\mathbb{C}^{2N}$ as
\[z_{01}=x_{01}+ix_{03}, z_{02}=x_{02}+ix_{04}, z_{k\ell}=x_{k\ell}+iy_{k\ell},\]
for $k=1,2,\ldots,\left\lfloor{N}/{2}\right\rfloor$ {and} $\ell=1,\ldots4.$

Based on our choice of coordinates and generators of $\D_N$ and the form of $e^{t{\cal M}}$, the action of $\D_N$ and $\T^m$ on the complex coordinates is
\begin{align*}
\begin{split}
\tilde \gamma \cdot z={}&\left(z_{01},z_{02},\exp\left(i\omega_1 \right)z_{11},\ldots,\exp\left(i\omega_1 \right)z_{14}, \ldots,\exp\left(i\omega_k \right)z_{k1},\ldots,\right.
\\
&\quad\left.\exp\left(i\omega_k\right)z_{k4}, \ldots, \exp\left(i\omega_{\left\lfloor {N}/{2}\right\rfloor}\right)z_{\left\lfloor{N}/{2}\right\rfloor 1},\ldots, \exp\left(i\omega_{\left\lfloor{N}/{2}\right\rfloor}\right)z_{\left\lfloor{N}/{2}\right\rfloor 4}\right),\\
\tilde\kappa \cdot z={}&\left(z_{01},z_{02},\overline z_{11},\ldots,\overline z_{14},\ldots,\overline z_{k1},\ldots,\overline z_{k4},\ldots,\overline z_{\left\lfloor {N}/{2} \right\rfloor 1},\ldots,\overline z_{\left\lfloor {N}/{2} \right\rfloor 4}\right),\AND\\
\tilde\theta \cdot z={}&\left(\exp\left(i\theta_0\right)z_{01},\exp\left(i\psi_0\right)z_{02},\exp\left(i\theta_1\right)z_{11},\exp\left(i\psi_1\right)z_{12},\exp\left(i\theta_1\right)z_{13},\right.\\
&\left.\exp\left(i\psi_1\right)z_{14},\ldots,\exp\left(i\theta_{\lfloor{N}/{2}\rfloor-1}\right)z_{(\lfloor{N}/{2}\rfloor-1) 1},\exp\left(i\psi_{\lfloor{N}/{2}\rfloor-1}\right)z_{(\lfloor{N}/{2}\rfloor-1) 2},\right.\\
&\exp\left( i\theta_{\lfloor{N}/{2}\rfloor -1}\right)z_{(\lfloor{N}/{2}\rfloor-1) 3}, \exp\left(i\psi_{\lfloor{N}/{2}\rfloor-1}\right)z_{(\lfloor{N}/{2}\rfloor-1)4},\ldots,\left. \exp\left(i\theta_{\lfloor{N}/{2}\rfloor}\right)z_{\left\lfloor {N}/{2}\right\rfloor 1},\right.\\
&\left. z_{\left\lfloor {N}/{2}\right\rfloor 2},\exp\left(i\theta_{\lfloor{N}/{2}\rfloor}\right)z_{\left\lfloor {N}/{2}\right\rfloor 3},z_{(\left\lfloor {N}/{2}\right\rfloor) 4}\right),
\end{split}
\end{align*}
where $\overline \cdot$ denotes the complex conjugate. Similarly, if $N$ is even, the generators are
\begin{equation*}
\tilde{\gamma}=P^{-1}\gamma P=\left(I_4, I_4\otimes \mathcal R(\omega_1),\ldots,I_4\otimes \mathcal R\left (\omega_{{N}/{2}-1}\right),-I_4\right)
\end{equation*}and
\begin{equation*}
\tilde{\kappa}=P^{-1}\kappa P=\left(I_4, I_4\otimes \mathcal S(\omega_1),\ldots,I_4\otimes \mathcal S\left (\omega_{{N}/{2}-1}\right),-I_4\right).
\end{equation*} In this case, we identify $X$
\[z_{01}=x_{01}+i x_{03}, z_{02}=x_{02}+ix_{04}, z_{k\ell}=x_{k\ell}+iy_{k\ell},\;\;\mbox{and}\;\;  z_{({N}/{2}) \ell}=x_{({N}/{2})\ell},\]
for $k=1,2,\ldots,{N}/{2}-1,$ and $\ell=1,\ldots4.$ In this case, the action on the complex coordinates becomes
\begin{align*}
\begin{split}
\tilde\gamma \cdot z={}&\left(z_{01},z_{02},\exp\left({i\omega_1}\right)z_{11},\ldots,\exp\left({i\omega_1}\right)z_{14}, \ldots,\exp\left({i\omega_{{N}/{2}-1}}\right)z_{({N}/{2}-1) 1},\ldots,\right.
\\
&\quad\exp\left({i\omega_{{N}/{2}-1}}\right)z_{({N}/{2}-1) 4},-z_{({N}/{2})1},\ldots,-z_{({N}/{2})4}\big),\\
\tilde\kappa \cdot z={}&\left(z_{01},z_{02},\overline z_{11},\ldots,\overline z_{14},\ldots,\overline z_{({N}/{2}-1)1},\ldots,\overline z_{({N}/{2}-1)4}, z_{({N}/{2})1},\ldots, z_{({N}/{2})4}\right),\AND
\\
\tilde\theta \cdot z={}&\left(\exp\left({ i\theta_0}\right)z_{01},\exp\left({ i\psi_0}\right)z_{02},\exp\left({ i\theta_1}\right)z_{11},\exp\left({ i\psi_1}\right)z_{12},\exp\left({ i\theta_1}\right)z_{13},\right.
\\
&\exp\left({i\psi_1}\right)z_{14},\ldots,\exp\left({ i\theta_{{N}/{2}-1}}\right)z_{({N}/{2}-1)1},\exp\left({ i\psi_{{N}/{2}-1}}\right)z_{({N}/{2}-1)2},
\\
&\quad\exp\left({ i\theta_{{N}/{2}-1}}\right)z_{({N}/{2}-1)3},\exp\left({ i\psi_{{N}/{2}-1}}\right)z_{({N}/{2}-1)4}, z_{({N}/{2})1},\ldots,z_{({N}/{2})4}\Big).
\end{split}
\end{align*}
Note that the Hamiltonian function in~\eqref{eq:hamiltonianNormalFormDegreeTwo} is already in normal form for all terms up to degree two as the linear normal forms commute with the $\D_N\times \Z_2\times \T^m$ actions above. We want to obtain the normal form up to degree four in~\eqref{eq:hamiltonianNormalFormDegreeTwo} and calculate the terms of degree four in $H_2$ which commute with $\D_N\times\Z_2\times \T^m$. This is done explicitly below.


 Let $u_{k\ell}=z_{k\ell}\overline z_{k\ell}$ and $v_{km}=z_{km}\overline z_{k(m+2)}$, then the degree two $\D_N$ invariants are
\begin{equation*}
u_{0\ell},\;\;u_{k\ell},\AND v_{km}+\overline v_{km},
\end{equation*}
for $k=1,\ldots, \left\lfloor{N}/{2}\right\rfloor$, $\ell=1,\ldots,4$, and $m=1,2$. For all $N$, the invariants $u_{0\ell},u_{k\ell}$ are also $\Z_2\times \T^m$ invariants. For $N$ odd,  we have $v_{\lfloor {N}/{2}\rfloor 2}+\ov{v_{\lfloor {N}/{2}\rfloor 2}}$ and for $N$ even $v_{({N}/{2}) m}+\ov{v_{({N}/{2}) m}}$ for $m=1,2$ as additional $\Z_2\times \T^m$ invariant.
The corresponding real invariants are
\begin{align}\label{eq:degree2InvariantsOdd}
\begin{split}
\mathcal U_1={}&x^2_{0m}+x^2_{0(m+2)},\\
\mathcal U_2={}&x^2_{k\ell}+y^2_{k\ell},\\
\mathcal U_3={}&x_{\lfloor {N}/{2}\rfloor 2}x_{\lfloor {N}/{2}\rfloor 4}+y_{\lfloor {N}/{2}\rfloor 2}y_{\lfloor {N}/{2}\rfloor 4}\qquad \mbox{$N$ odd},\\
\mathcal U_4={}&x_{(N/2)l}x_{(N/2)l}\qquad \mbox{$N$ even},
\end{split}
\end{align}
where $m=1,2$, $k=1,\ldots,\lfloor{N}/{2}\rfloor$, $\ell=1,2,3,4$. Degree four invariants are calculated as the products of the degree two invariants. The list  of possible invariants of the system is long, the ones relevant to the gyroscopic system can be found in~\ref{sec:degree4invariants}.
We represent the Hamiltonian function in normal form up truncated to degree four as
\begin{equation}\label{eq:hamiltonianNormalFormDegreeFour}
\tilde H(X)=\tilde H_0(X)+\tilde H_2(X).  
\end{equation}
The local dynamics near the bifurcation point can be studied via this normal form. However, because of the large number of terms appearing in $\tilde H_2(X)$, this is a cumbersome exercise. We obtain the explicit normal in the case study $N=3$ in a section below. However, for the general case, we determine the nature of the bifurcation using the splitting lemma which preserves the zero set of the vector field, but not the local dynamics.

\subsection{Splitting Lemma}
With the Hamiltonian function now in normal form, we can further simplify the system by applying the equivariant splitting lemma~\cite{bridges1993singularity}. This simplification allows us to separate the degenerate and nondegenerate variables of the Hamiltonian and thus find the essential nonlinear terms necessary for further analysis.

Let $f:\mathbb{R}^n\to\mathbb{R}$ be a  $\Gamma$-equivariant function. A critical point of $f$ is not degenerate if the determinant of its Hessian matrix is nonzero and it is degenerate otherwise. Suppose $x_0=0$ is a degenerate singular point of $f$ and the corresponding Hessian matrix has rank of $m$ and corank of $k$, where $n=m+k$. For a function with the aforementioned qualities, the equivariant splitting lemma states that there must exist a change of coordinates in the neighborhood of the critical point such that
$$f(x(\chi, u),u)=K(\chi)+h(u),$$
where $\chi\in\mathbb{R}^m$, $u\in\mathbb{R}^k$, $K$ is the restriction of $\displaystyle\frac{1}{2}d^2f$ to $\mathbb{R}^m\times \{0\}$, and $h$ is the  remainder function. This remainder function $h$ can be found implicitly. For each $u$ near the origin,  there is a unique point $x=\chi(u)$ such that $d_{X}\tilde H( \chi(u),u)=0$ and
$$h(u)=\tilde H(\chi(u),u).$$ Thus, we solve $d_{X}\tilde H(\chi(u),u)=0$ for each component of $\chi$ and substitute each $\chi_i$ in terms of $u_i$ back into $\tilde H(\chi(u),u)$.

Considering the Hamiltonian function in~\eqref{eq:hamiltonianNormalFormDegreeFour}, we may write it as
$$\tilde H(X)=\tilde H_0(X) +\tilde H_2(X)$$  by ignoring the higher order terms. There is a degenerate critical point for $\tilde H(X)$ at $X=0$. Based on the normal form of $H_0$ in~\eqref{eq:degreeTwoHamiltonianNormalFormOdd}, when $N$ is odd, the rank of the corresponding Hessian matrix is $N-2$, and the corank is two. By the equivariant splitting lemma, there exist a change of coordinates in a neighborhood of the origin such that
$$ \tilde H(x(\chi,u),u)=K(\chi)+ h(u),$$
where $\chi \in\mathbb{R}^{N-2}$, $u\in\mathbb{R}^2$. Examining the Hessian matrix, we find that the restriction of $\tilde H(X)$ to $\mathbb{R}^{N-2}\times\left\{(x_{\lfloor N/2\rfloor 4},y_{\lfloor N/2\rfloor 4})=(0,0)\right\}$ has a
nondegenerate critical point at $X_0=0$. Thus, we may write $\chi=(\chi_1,\ldots,\chi_{N-2})$ and $u=(u_1,u_2)$ as
\begin{align*}
\chi={}&(x_{01},\ldots,x_{04}, \ldots, y_{(\left\lfloor {N}/{2}\right\rfloor -1)1},\ldots,\\
& \quad y_{(\left\lfloor {N}/{2}\right\rfloor -1)4},x_{\left\lfloor {N}/{2}\right\rfloor 1},x_{\left\lfloor {N}/{2}\right\rfloor 2},x_{\left\lfloor {N}/{2}\right\rfloor 3},y_{\left\lfloor {N}/{2}\right\rfloor 1},y_{\left\lfloor {N}/{2}\right\rfloor 2},x_{\left\lfloor {N}/{2}\right\rfloor 3}),
\qquad \AND \\
 u={}&(x_{\lfloor N/2\rfloor 4},y_{\lfloor N/2\rfloor 4}).
 \end{align*}
A direct calculation shows that $\chi=0$ is always a solution to $d_{X}\tilde H(\chi(u),u)=0$. Thus the remainder function is
\begin{equation*}\label{eq:reminderFunctionOdd}
h^{odd}(u)= h^{odd}(x_{\left\lfloor {N}/{2}\right\rfloor 4},y_{\left\lfloor {N}/{2}\right\rfloor 4})=\alpha_1\left(x_{\left\lfloor {N}/{2}\right\rfloor 4}^{2}+x_{\left\lfloor {N}/{2}\right\rfloor 4}^{2}\right)^2,
\end{equation*}where $\alpha_1$ is a constant in terms of $\mu,\kappa$, and $\Omega$.

When $N$ is even, the corresponding Hessian has rank of $N-1$ and corank of one. Thus, following the same steps for $N$ is odd, the remainder function is
\begin{equation*}\label{eq:reminderFunctionEven}
h^{even}(u)= h^{even}(y_{\left\lfloor {N}/{2}\right\rfloor 4})=\alpha_2 y_{\left\lfloor {N}/{2}\right\rfloor 4}^{4},
\end{equation*}where $\alpha_2$ is another constant in terms of $\mu,\kappa$, and $\Omega$.

\section{Case Study: $\D_3$-Symmetric Gyroscopic System}
\label{sec:bifurcationParameter}
In Sections~\ref{sec:hamiltonianFormulation} through~\ref{sec:nonlinearNormalForm} we
have analyzed the collective behavior of $N$ gyroscopes bidirectionally coupled in a ring fashion
through their driving axes. In this section, we illustrate the general theory by studying, in particular,
a relatively small ring consisting of $N=3$ gyroscopes, so that the system exhibits $\D_3$ symmetry.

\subsection{The $\D_3$-Symmetric System}
\label{sec:theD3System}
Assuming that we have performed the isotypic decomposition outlined in Section~\ref{sec:isotypicDecomposition} and using the same notations as before, the $D_3$ symmetric system can be written as
\begin{equation*}
\dot U = \mathcal M U+ F,
\end{equation*} where $U=(U_1,U_2,U_3)^T$, $\mathcal M=\operatorname{diag} (M_1+2M_2,M_1-M_2,M_1-M_2)$ and $F=(F_1,F_2,F_3)^T$.

Clearly,  the eigenvalues of the system are the eigenvalues of $M_1+2M_2$ and $M_1-M_2$.  The four roots corresponding to the characteristic polynomial  of the $M_1+2M_2$ block are $\pm\sqrt{-2\Omega^2-\kappa\pm 2\Omega\sqrt{\Omega^2+\kappa}} $.
Since $2\Omega^2+\kappa- 2\Omega\sqrt{\Omega^2+\kappa}$ is greater than zero if and only if $\kappa^2>0$, there are always two pairs of purely imaginary eigenvalues for the $M_1+2M_2$ block and they are
\begin{equation}\label{eq:eigenvaluesForD3System}
\pm i\sqrt{2\Omega^2+\kappa+2\Omega\sqrt{\Omega^2+\kappa}}\AND\pm i\sqrt{2\Omega^2+\kappa-2\Omega\sqrt{\Omega^2+\kappa}}.
\end{equation}

For the $M_1-M_2$ block, the roots of the corresponding characteristic polynomial  are
\[
\pm\sqrt{-\kappa-2\,{\Omega}^{2}-\frac{3}{2}\,\lambda\pm\frac{1}{2}\,\sqrt {16\,{\Omega}^{2}\kappa+16\,{\Omega}^{4}+24\,{\Omega}^{2}\lambda+9\,{\lambda}^{2}}
}\, .
\]
Since $-\kappa-2\,{\Omega}^{2}-3/2\,\lambda-1/2\,\sqrt {16\,{\Omega}^{2}\kappa+16\,{\Omega}^{4}+24\,{\Omega}^{2}\lambda+9\,{\lambda}^{2}}
$ is real and negative for positive parameter values, one set of the eigenvalues must be a purely imaginary pair of the form
\[
\pm i\sqrt{\kappa+2\,{\Omega}^{2}+\frac{3}{2}\,\lambda+\frac{1}{2}\,\sqrt {16\,{\Omega}^{2}\kappa+16\,{\Omega}^{4}+24\,{\Omega}^{2}\lambda+9\,{\lambda}^{2}}
}\, ,
\]
and the other set of eigenvalues are
\[
\pm\sqrt{-\kappa-2\,{\Omega}^{2}-\frac{3}{2}\,\lambda+\frac{1}{2}\,\sqrt {16\,{\Omega}^{2}\kappa+16\,{\Omega}^{4}+24\,{\Omega}^{2}\lambda+9\,{\lambda}^{2}}
}\, .
\]
Setting $-\kappa-2\,{\Omega}^{2}-3/2\,\lambda+1/2\,\sqrt {16\,{\Omega}^{2}\kappa+16\,{\Omega}^{4}+24\,{\Omega}^{2}\lambda+9\,{\lambda}^{2}}=0$, this expression simplifies to
$\kappa(\kappa+3\lambda)=0$. If $\lambda>-\ds {\kappa \over 3}$ then the pair of eigenvalues is purely
imaginary and it switches to a pair of real eigenvalues with opposite sign as $\lambda$ crosses the
critical value $\lambda^*=-\ds {\kappa \over 3}$. We wish to point out that this is the same critical value
of the coupling strength that was found via perturbation analysis in~\cite{vu2010two}.

\subsection{Symplectic Transition Matrices}
By directly applying the results from Section~\ref{sec:linearAnalysis}, the components of the diagonal
symplectic matrix $Q=\operatorname{diag}(Q_0,Q_1,Q_1)$ are found to be
\[
Q_0=
\left[
\begin {array}{cccc}
\sqrt{\frac{\nu_1}{2\xi_1}}&0&0&\sqrt{\frac{\nu_2}{2\xi_2}}
\\
0&\frac{2\Omega+\sqrt{\Omega^2+\kappa}}{\sqrt{2\xi_2\nu_2}}&\frac{2\Omega-\sqrt{\Omega^2+\kappa}}{\sqrt{2\xi_1\nu_1}}&0
\\
0&\sqrt{\frac{\xi_2}{2\nu_2}}&\sqrt{\frac{\xi_1}{2\nu_1}}&0
\\
\sqrt{\frac{(4\Omega^2+\kappa)\nu_1}{2\xi_1}}&0&0&\sqrt{\frac{(4\Omega^2+\kappa)\nu_2}{2\xi_2}}
\end {array}
\right]
\]
and
\[
Q_1= \left[ \begin {array}{cccc}
0&0&\,-{\frac {4\Omega}{ \left( \kappa+16\,{\Omega}^{2} \right) ^{3/4}}}&{\frac {\kappa}{\sqrt {\kappa\, \left( \kappa+16\,{\Omega}^{2} \right) }}}
\\
{\frac {1}{\sqrt [4]{\kappa+16\,{\Omega}^{2}}}}&\,{\frac {4\Omega}{\sqrt {\kappa\, \left( \kappa+16\,{\Omega}^{2} \right) }}}&0&0\\
{\frac {2\Omega}{\sqrt [4]{\kappa+16\,{\Omega}^{2}}}}&-{\frac {\kappa+8\,{\Omega}^{2}}{\sqrt {\kappa\, \left( \kappa+16\,{\Omega}^{2} \right) }}}&0&0\\
0&0&{\frac {\kappa+8\,{\Omega}^{2}}{ \left( \kappa+16\,{\Omega}^{2} \right) ^{3/4}}}&\,{\frac {2\kappa\,\Omega}{\sqrt {\kappa\, \left( \kappa+16\,{\Omega}^{2} \right) }}}\end {array} \right]
\]
where $\nu_1=\sqrt{2\Omega^2+\kappa-2\Omega\sqrt{\kappa+\Omega^2}}$, $\nu_2=\sqrt{2\Omega^2+\kappa+2\Omega\sqrt{\kappa+\Omega^2}}$, $\xi_1=\Omega^2+\kappa-\Omega\sqrt{\Omega^2+\kappa}$, and $\xi_2=\Omega^2+\kappa+\Omega\sqrt{\Omega^2+\kappa}$.
Using the symplectic transformation $Q$, the linear part of system becomes $\mathsf M=\operatorname{diag}(\mathsf M_0,\mathsf M_1,\mathsf M_1)$, where
\[
\mathsf M_0= Q_0^{-1}(M_1+2M_2)Q_1=
\left(\begin{array}{cccc}
0&0&\nu_1&0\\
0&0&0&\nu_2\\
-\nu_1&0&0&0\\
0&-\nu_2&0&0
\end{array}\right)
\]
and
 \[
 \mathsf M_1= Q_1^{-1}(M_1-M_2)Q_1=
\left(\begin{array}{cccc}
0&0&\sqrt{\kappa+4\Omega^2}&0\\
0&0&0&0\\
-\sqrt{{\kappa+4\Omega^2}}&0&0&0\\
0&-1&0&0
\end{array}\right).
\]

\subsection{Hamiltonian Function}
With the linear part of the system in normal form, we proceed to put the high order terms in normal form as well. Suppose $X=(X_0,X_1,Y_1)\in\mathbb{R}^{12}$, then let $U= QX$. The Hamiltonian function $H$ can now be written as
\begin{equation*}\label{eq:hamiltonianNormalFormDegreeTwoD3}
H(X)=\tilde H_0(X)+H_2(X),
\end{equation*}
where $\tilde H_0(X)$ and $H_2(X)$ represent  polynomials of degree two and four, respectively, Furthermore, $\tilde H_0(X)$ is already in normal and it is
\begin{equation*}
\begin{split}
\tilde H_0(X)={}& \frac{1}{2}X^TJ^{-1}AX\\
={}&\frac{\nu_1}{2}(x_{01}^2+x_{03}^2)+\frac{\nu_2}{2}(x_{02}^2+x_{04}^2)+\frac{1}{2}(x_{12}^2+y_{12}^2)+\frac{\sqrt{\kappa+4\Omega^2}}{2}(x_{11}^2+y_{11}^2)
\\
{}&+\frac{\sqrt{\kappa+4\Omega^2}}{2}(x_{13}^2+y_{13}^2).
\end{split}
\end{equation*}
The expression for $H_2$ is too long to be reproduced in full. Based on the results in Section~\ref{sec:nonlinearNormalForm} and~\ref{sec:degree4invariants}, we know that
\begin{align*}
\begin{split}
g_1&=\left(x_{11}^2+y_{11}^2\right)^2,\\
g_2&=\left(x_{11}^2+y_{11}^2\right)\left(x_{12}^2+y_{12}^2\right),\\
g_3&=\left(x_{12}^2+y_{12}^2\right)^2,\\
g_4&=\left(x_{13}^2+y_{13}^2\right)^2,\\
g_5&=\left(x_{13}^2+y_{13}^2\right)\left(x_{14}^2+y_{14}^2\right),\;\mbox{and}\\
g_6&=\left(x_{14}^2+y_{14}^2\right)^2,\\
 \end{split}
 \end{align*}
 are the only relevant invariants that occur in $H_2(X)$. Thus, the normal form for $H_2(X)$ can be written as
 $$\tilde H_2(X)=\mu\left(a_1g_1+a_2g_2+a_3g_3+a_4g_4+a_5g_5+a_6g_6\right),$$
where
$\ds a_1=\frac{1}{8}\,{\frac {1}{\kappa+4\,{\Omega}^{2}}}$,
$\ds a_2=\,{\frac {{\Omega}^{2}}{\kappa \left( \kappa+4\,{\Omega}^{2} \right) ^{3/2}}}$,
$\ds a_3=\,{\frac {{2\,\Omega}^{4}}{\kappa^2\left( \kappa+4\,{\Omega}^{2} \right) ^{2}}}$,
$\ds a_4=\,{\frac {2\,{\Omega}^{4}}{ \left( \kappa+4\,{\Omega}^{2} \right) ^{3}}}$,
$\ds a_5=\,{\frac {\kappa\,{\Omega}^{2}}{ \left( \kappa+4\,{\Omega}^{2} \right) ^{5/2}}},$ and
$\ds a_6=\frac{1}{8}\,{\frac {\,{\kappa}^{2}}{ \left( \kappa+4\,{\Omega}^{2} \right) ^{2}}}$.
\vspace*{0.1in}

Setting $\epsilon=\mu$, the Hamiltonian function in normal form is
\begin{equation}\label{eq:hamiltonianFucntionNormalForm}
   \tilde H(X)=\tilde H_0(X)+\tilde H_2(X)+\mathcal O(\epsilon^2).
\end{equation}

We now directly apply the equivariant splitting lemma to further simplify the Hamiltonian function. Suppose $\chi=(\chi_1,\ldots,\chi_{10})$ and $u=(u_1,u_2)$, then $$\chi=(x_{01},\ldots,x_{04},x_{11},x_{12},x_{13},y_{11},y_{12},y_{13})\AND u=(x_{14},x_{14}).$$ We solve $d_{X}\tilde H(\chi(u),u)=0$ for each $\chi_i$ in terms of $u_i$. Back substituting into $\tilde H(\chi(u),u)$, we find that
$$h(u)= h(x_{14},y_{14})=\frac{1}{8}\,{\frac {\mu\,{\kappa}^{2}}{ \left( \kappa+4\,{\Omega}^{2} \right) ^{2}}}\left(x_{14}^{2}+y_{14}^{2}\right)^2.$$

\subsection{Introducing a Bifurcation Parameter}
Note that the normal form obtained in Section~\ref{sec:linearAnalysis} was calculated at the critical value of the coupling strength. Thus the normal form reduction is primarily valid at criticality but it cannot provide information on the system dynamics away from the critical point.  To overcome this deficiency, we will introduce a bifurcation parameter to study the dynamics in a neighborhood of the critical point.

As mentioned in Section~\ref{sec:linearAnalysis}, the critical point occurs at $\ds \lambda_c=-\frac{1}{3}\kappa$. Let $|\tilde \eta|\ll0$ and add it to the critical coupling strength, so that  $\lambda_c=-\frac{1}{3}\kappa+\tilde\eta$. For notational convenience, we rescale the new parameter as $\ds \tilde\eta=\frac{1}{3}\eta$. Let
\[
\Upsilon=\left(\begin{array}{cccc}
0&0&0&0\\
0&0&0&0\\
-1&0&0&0\\
0&0&0&0
\end{array}\right),
\]
then the linear perturbation of the system can be written as $M+\eta \Phi$, where
$$
\Phi=\left(\begin{array}{ccc}
0_4&0_4&0_4\\
0_4&\Upsilon&0_4\\
0_4&0_4&\Upsilon
\end{array}\right)
$$ with $0_4$ as the $4\times4$ zero matrix.

A direct calculation shows that $\Phi$ is a Hamiltonian matrix. Suppose matrices $P$ and $Q$ are as described in Sections~\ref{sec:isotypicDecomposition} and~\ref{sec:linearAnalysis}. Then we may write the Hamiltonian function associated
with $\Phi$ as
\begin{align}\label{eq:perturbedHamiltonian}
\begin{split}
\eta H_0^P&=\frac{\eta}{2}X^TJ^{-1}(QP)^{-1}N(QP)X\\
&=b_1x_{01}^2+b_2x_{04}^2+b_3x_{13}^2+b_4x_{14}^2+b_5y_{13}^2+b_6y_{14}^2,
\end{split}
\end{align}where
$\ds b_1=\frac{1}{6}\,{\frac { \xi_2 \nu_1}{\kappa\, \left( \kappa+\,{\Omega}^{2} \right) }}$,
$\ds b_2=\frac{1}{6}\,{\frac { \xi_1 \nu_2}{\kappa\, \left( \kappa+\,{\Omega}^{2} \right) }}$,
$\ds b_3=\,{\frac {{\Omega}^{2}}{ \left( \kappa+4\,{\Omega}^{2} \right) ^{3/2}}}$,
$\ds b_4=\frac{1}{4}\,{\frac {\kappa}{\kappa+4\,{\Omega}^{2}}}$,
$\ds b_5={\frac {5}{3}}\,{\frac {{\Omega}^{2}}{ \left( \kappa+4\,{\Omega}^{2} \right) ^{3/2}}}$, and
$\ds b_6={\frac {5}{12}}\,{\frac {\kappa}{\kappa+4\,{\Omega}^{2}}}$.\\

Comparing the  terms in~\eqref{eq:perturbedHamiltonian} to the list of degree two invariants in~\eqref{eq:degree2InvariantsOdd}, the Hamiltonian function of the linear perturbation in normal form is
 $$\tilde {H}_0^P=b_3\left( x_{13}^{2}+y_{13}^{2} \right)+b_4\left( x_{14}^{2}+y_{14}^{2} \right).$$ Thus the Hamiltonian function corresponding to the linear part of the differential system is
\begin{equation*}
H_L=\tilde H_0+\eta\tilde{H}_0^P.
\end{equation*}
One may derive the corresponding Jacobian from the linear Hamiltonian function. The eigenvalues
of the Jacobian pertaining to $H_L$ are
\begin{equation*}
  \begin{array}{lll}
   \lambda_1 &=& \ds \pm\nu_1, \\
   \lambda_2 &=& \ds \pm\nu_2, \\[10pt]
   \lambda_3 &=& \ds {\frac {\sqrt {- \left( \kappa+4\,{\Omega}^{2} \right)  \left( {\kappa}^{2}+8\,{\Omega}^{2}\kappa+16\,{\Omega}^{4}+2\,\eta\,{\Omega}^{2} \right) }}{\kappa+4\,{\Omega}^{2}}},\AND \\
   \lambda_4 &=& \ds {\frac {\sqrt {- \left( 2\,\kappa+8\,{\Omega}^{2} \right) \eta\,\kappa}}{\kappa+4\,{\Omega}^{2}}},
   \end{array}
\end{equation*}
where $\lambda_3$ and $\lambda_4$ both have algebraic multiplicity of four. Clearly, $\lambda_1$ and
$\lambda_2$ are the same eigenvalues for the $M_1+2M_2$ block found in~\eqref{eq:eigenvaluesForD3System} and they are unaffected by values of $\eta$. Furthermore, we observe that at the critical value of
$\eta=0$, these eigenvalues are the same as the eigenvalues found in Section~\ref{sec:theD3System}. Regardless of the value of $\eta$, $\lambda_3$ must be purely imaginary. When $\eta<0$, $\lambda_4$ must be real. As $\eta$ increases and becomes zero, the eigenvalues also become zero. After $\eta$ crosses criticality and becomes positive,
$\lambda_4$ becomes purely imaginary.

\subsection{Numerical Simulations}


Based on the results from the splitting lemma, and after adding the perturbation term, we may restrict the domain of the Hamiltonian function to $x_{13},x_{14}, y_{13}$ and $y_{14}$ and the reduced Hamiltonian is
\begin{equation} \label{eq:reduced_H_function}
\mathfrak H=\tilde H_0|_{(x_{13},x_{14},y_{13},y_{14})}+\eta\tilde H_0^P+h(x_{14},y_{14}).
\end{equation}
Thus, the equations of motion are
\begin{align} \label{eq:reduced_Hamiltonian}
\begin{split}
\frac{\partial\mathfrak H}{\partial x_{13}}&=x_{13}\sqrt {\kappa+4\,{\Omega}^{2}}+2\,{\frac {\eta\, {\Omega}^{2}x_{13}}{ \left( \kappa+4\,{\Omega}^{2} \right) ^{3/2}}},
\\
-\frac{\partial\mathfrak H}{\partial x_{14}}&=\frac{1}{2}\,{\frac {\eta\,\kappa\,x_{14}}{\kappa+4\,{\Omega}^{2}}}+\frac{1}{2}\,{\frac {\mu\,{\kappa}^{2}x_{14}^{3}}{ \left( \kappa+4\,{\Omega}^{2} \right) ^{2}}}
+\frac{1}{2}\,{\frac {\mu\,{\kappa}^{2}x_{14}y_{14}^{2}}{ \left( \kappa+4\,{\Omega}^{2} \right) ^{2}}},
\\
\frac{\partial\mathfrak H}{\partial y_{13}}&=y_{13}\sqrt {\kappa+4\,{\Omega}^{2}}+2\,{\frac {\eta\, {\Omega}^{2} y_{13}}{ \left( \kappa+4\,{\Omega}^{2} \right) ^{3/2}}},
\\
-\frac{\partial\mathfrak H}{\partial y_{14}}&=\frac{1}{2}\,{\frac {\eta\, \kappa y_{14}}{\kappa+4\,{\Omega}^{2}}}+\frac{1}{2}\,{\frac {\mu\,{\kappa}^{2}y_{14}^{3}}{ \left( \kappa+4\,{\Omega}^{2} \right) ^{2}}}
+\frac{1}{2}\,{\frac {\mu\,{\kappa}^{2}x_{14}^{2}y_{14}}{ \left( \kappa+4\,{\Omega}^{2} \right) ^{2}}}.
\end{split}
\end{align}

\begin{figure}
   \begin{center}\begin{subfigure}[b]{0.49\textwidth}
\centering
\includegraphics[width=\textwidth]{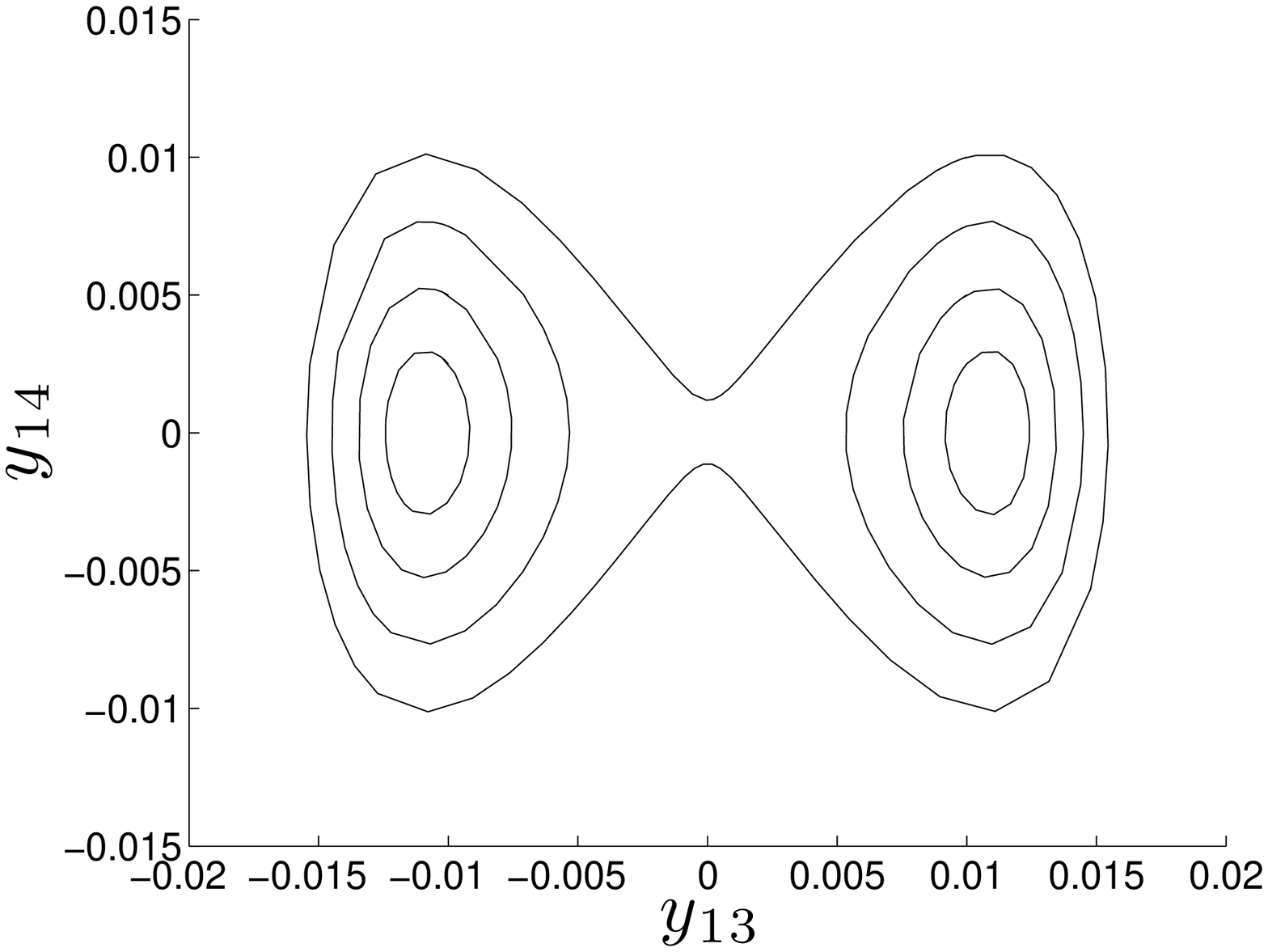}
\caption{}
\label{fig:negative}
\end{subfigure}
\begin{subfigure}[b]{0.49\textwidth}
\centering
\includegraphics[width=\textwidth]{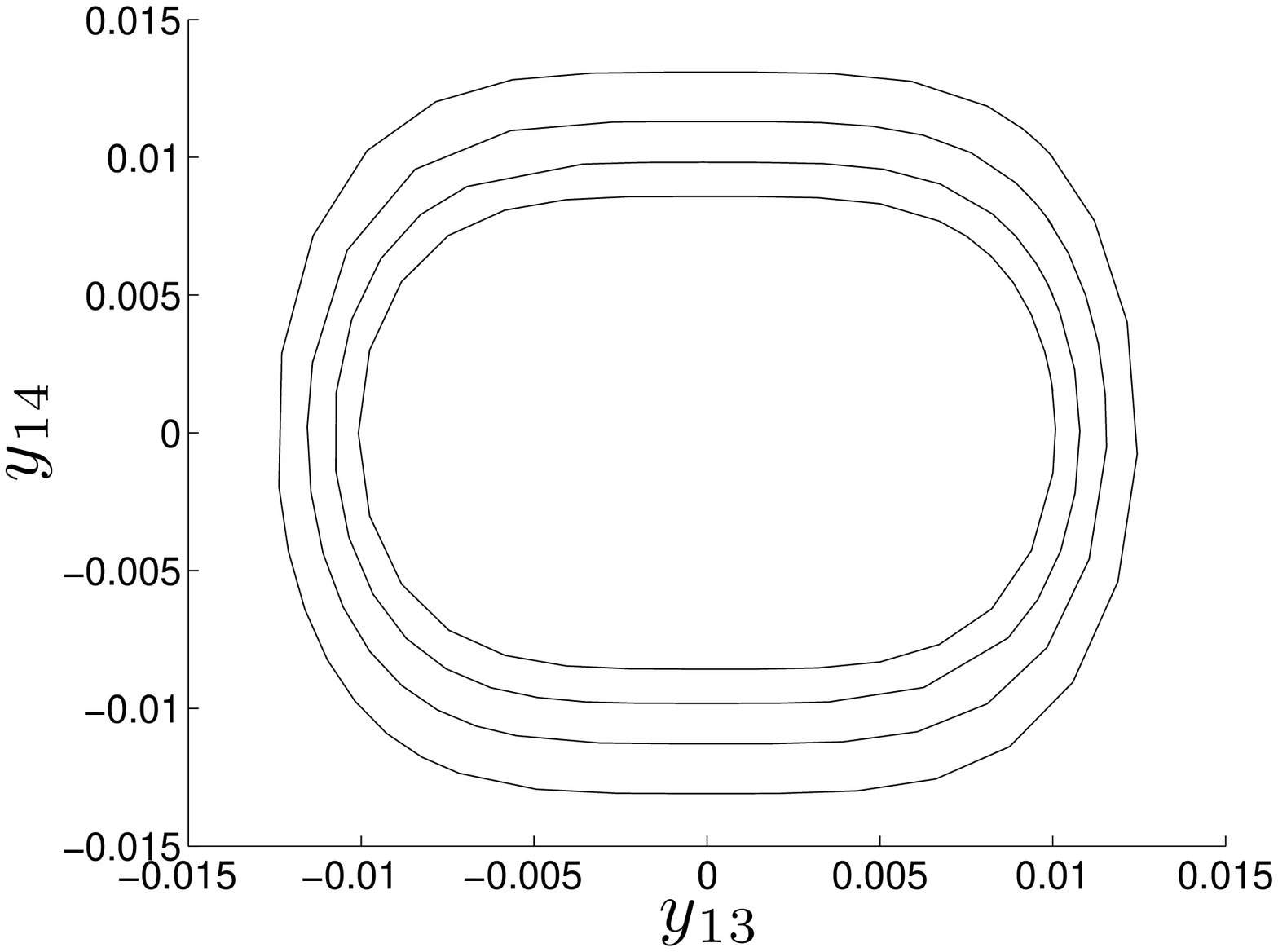}
\caption{}
\label{fig:positive}
\end{subfigure}
 \end{center}
  \caption{\small Computer simulations of the reduced Hamiltonian system~\Ref{eq:reduced_Hamiltonian}
     successfully capture the transition in the dynamics of a $\D_3$-symmetric coupled gyroscope
     system as the coupling strength varies. In Figure~\ref{fig:negative}, when $\eta<0$, a pair of stable centers arise and each
     one is surrounded by a family of periodic oscillations. In the full system~\Ref{eq:coupled_gyro},
     two of the driving modes are completely synchronized. They oscillate with either a positive or negative
     mean, which correspond to the positive or negative values of the centers in the figure. The third mode
     oscillates in phase with respect to the other two but with the opposite sign in the mean oscillations.
     As $\eta$ increases, the two centers move closer towards the saddle-point until
    they all coalesce at $\eta=0$. We observe in Figure~\ref{fig:positive} that when $\eta>0$, there is only one stable center and it is     surrounded by a family of stable periodic oscillations. In the full system, the non-zero mean  oscillations no longer exist and only one stable periodic oscillation is observed, i.e., the
     complete synchronization has now emerged via a pitchfork bifurcation. Parameters are as in table~\ref{tab:parameters} with $\Omega=308$.}
   \label{fig:simulations}
\end{figure}

Computer simulations of the coupled gyroscope dynamics, as is captured by the reduced Hamiltonian
system in~\Ref{eq:reduced_Hamiltonian}, were carried out with parameter values assigned according to
table~\ref{tab:parameters}. Based on a comparison of the reduced Hamiltonian
function in~\Ref{eq:reduced_H_function} to the bifurcations types listed in~\cite{buono2005symmetric}, it
is reasonable to expect a super critical pitchfork bifurcation to occur in the coupled gyroscope dynamics
as $\eta$ crosses criticality, i.e, $\eta=0$, which can also be interpreted as
$\lambda=\lambda_c$. This is indeed the case and the actual transition is illustrated in
figure~\ref{fig:simulations}. When $\eta < 0$, the phase space dynamics exhibits, see
figure~\ref{fig:simulations}(left), a pair of stable centers (one positive and one negative) each one
surrounded by a family of periodic oscillations and an unstable saddle point at zero. As $\eta$
increases, the centers get closer to one another and to the saddle-point at zero until, eventually, at
$\eta=0$ they all coalesce into a single center still surrounded by a family of periodic solutions, as
is shown in~figure~\ref{fig:simulations}(right). This entire transition corresponds to the pitchfork bifurcation
that leads to complete synchronization in the original coordinates of the full system~\Ref{eq:coupled_gyro},
as it was reported in~\cite{vu2010two}. That is, when $\lambda < \lambda_c$ there are two types of
periodic patterns. One unstable synchronized state, in which all driving-mode oscillations are in phase
and they all oscillate with the same amplitude and zero-mean. And one stable pattern where two of the
driving modes are completely synchronized, oscillating with either a positive or negative mean (which
corresponds to the positive/negative centers), while the third mode oscillates in phase with respect to
the other two but with the opposite sign in the mean of the oscillations. This stable pattern can also be
described as two gyroscopes oscillating around one of the two wells of the energy function represented
by the Hamiltonian function in~\Ref{eq:reduced_H_function} and the third one oscillating around the other
well. As $\lambda$ approaches $\lambda_c$ the absolute value of the mean oscillations of the stable
pattern gradually decreases until it becomes zero at $\lambda = \lambda_c$. Passed $\lambda_c$, the
non-zero mean oscillations disappear while the complete synchronization state with zero-men
oscillations becomes locally asymptotically stable.

\section{Discussion and Conclusion}
 \label{sec:conclusion}
Ideas and methods from equivariant bifurcation theory were used to study the equations of motion of a
high-dimensional coupled nonlinear system with Hamiltonian structure. The equations belong to a
particular model for a gyroscope system but the theory developed in this work is generic enough to
study a wider range of coupled Hamiltonian systems with symmetry.
Coupling among the individual systems lead to high dimensionality and, in some cases, the specific
choice of coupling function can destroy the Hamiltonian structure. For instance, a ring array with
nearest-neighbor coupling with a preferred orientation, i.e., unidirectional coupling, leads to a network
with global $\Z_N$-symmetry, where $\Z_N$ is the group of cyclic rotations of $N$ objects. If there is
no preferred orientation, i.e., bidirectional coupling, then the ring possesses $\D_N$ symmetry, where
$\D_N$ is the dihedral group of symmetries of a regular $N$-gon. It was found that in the former case,
the $\Z_N$-symmetry actually destroys the Hamiltonian structure while in the latter, the $\D_N$-symmetry
preserves the Hamiltonian structure. An interesting question that arises almost immediately is to
determine the type of coupling functions that can preserve the Hamiltonian structure for a generic
network of coupled nonlinear systems, e.g., nonlinear oscillators. A complete answer to this question
should include linear as well as nonlinear coupling functions and the task is referred for future work.
Symplectic transformations were calculated to rewrite the linear and nonlinear terms of the network
equations in normal form and to facilitate a bifurcation analyses of the network equations valid for
any ring size $N$. The analysis produced an analytical expression for the critical value of the
coupling strength that leads to completely synchronized behavior, i.e., same amplitude and phase
of oscillations, that is also valid for any ring size. This result is significant because synchronization
leads to improved performance and robustness against phase-drift and, thus, knowledge of the
critical coupling parameter is important to aid in the design and operation of an actual device. The
results of the generic theory were then illustrated with a particular ring with $\D_3$ symmetry. In this
case, the reduced Hamiltonian function, via normal forms, successfully capture the nature of the
pitchfork bifurcation that leads the three-ring system to synchronize as it was previously reported
through perturbation analysis. It is our hope that the analysis presented in this manuscript can
lead to a better understanding of the role of symmetry in many other highly-dimensional generic
systems with Hamiltonian structure.

\section*{Acknowledgments}
B.S.C. and A.P. were supported by the Complex Dynamics and Systems Program of the Army
Research Office, supervised by Dr. Samuel Stanton, under grant W911NF-07-R-003-4.
A.P. was also supported by the ONR Summer Faculty Research Program, at SPAWAR Systems
Center, San Diego. V.I. acknowledges support from the Office of Naval Research (Code 30) and
the SPAWAR internal research funding (S\&T) program.
We also would like to acknowledge constructive discussions with Dr. Brian Meadows at SPAWAR
Systems Center Pacific and with Prof. Takashi Hikihara and Dr. Suketu Naik at Kyoto University.
P-L.B. would like to thank Alberto Alinas for checking some early calculations as part of a student
project. P-L.B. acknowledges the funding support from NSERC (Canada) in the form of a Discovery
Grant.

\newpage
\appendix
\section{Calculation of Invariants}
\label{sec:degree4invariants}
There are many possible higher order invariants based on the symmetry of the system, but only some of them are relevant to the study of the coupled gyroscopic system. Those are extracted in this section.
Since we have already detailed the calculations of the degree two invariants and degree three terms do not appear in the Hamiltonian function,  we consider the degree four terms in $H_2$.
Based on the symplectic matrices $P$ and $Q$ found in Sections~\ref{sec:isotypicDecomposition} and~\ref{sec:linearAnalysis}, we may write the coordinate transformation as
$$\left(\begin{array}{c}q_1\\p_1\\\vdots\\q_N\\p_N\end{array}\right)=QP\left(\begin{array}{c}X_1\\\vdots\\\vdots\\\vdots\\X_N\end{array}\right).$$
Thus, we may think of the  coordinate and momentum variables $q_i$ and $p_i$ as functions of $X_i=\left(x_{i,1},x_{i,2},x_{i,3},x_{i,4}\right)$.

When $N$ is odd,  the matrix $P$ may be written in block matrix form as
$$P=\left(\begin{array}{ccccccccc}
I_4&I_4&I_4&\ldots&I_4&I_4\\
I_4&\Im(\zeta)I_4&\Re\left(\zeta\right)I_4&\ldots&\Im\left(\zeta^{N-1}\right)I_4&\Re\left(\zeta^{N-1}\right)I_4\\
\vdots&\vdots&\vdots&&\vdots&\vdots\\
I_4&\Im\left(\zeta^{N-1}\right)I_4&\Re\left(\zeta^{N-1}\right)I_4&\ldots&\Im\left(\zeta^{(N-1)\left\lfloor {N}/{2}\right\rfloor}\right)I_4&\Re\left(\zeta^{(N-1)\left\lfloor {N}/{2}\right\rfloor}\right)I_4
\end{array}\right)$$
and $Q$ may also be written as $Q=\operatorname {diag}\left(Q_0,Q_1,Q_1,\ldots,Q_{\left\lfloor{N}/{2}\right\rfloor},Q_{\left\lfloor{N}/{2}\right\rfloor}\right)$. The product of these two matrices is
$$
\footnotesize
QP=\left(\begin{array}{ccccccccc}
Q_0&Q_1&Q_1&\ldots&Q_{\left\lfloor {N}/{2}\right\rfloor}&Q_{\left\lfloor {N}/{2}\right\rfloor}\\
Q_0&\Im(\zeta)Q_1&\Re\left(\zeta\right)Q_1&\ldots&\Im\left(\zeta^{N-1}\right)Q_{\left\lfloor {N}/{2}\right\rfloor}&\Re\left(\zeta^{N-1}\right)Q_{\left\lfloor {N}/{2}\right\rfloor}\\
Q_0&\Im\left(\zeta^2\right)Q_1&\Re\left(\zeta^2\right)Q_1&\ldots&\Im\left(\zeta^{2(N-1)}\right)Q_{\left\lfloor {N}/{2}\right\rfloor}&\Re\left(\zeta^{2(N-1)}\right)Q_{\left\lfloor {N}/{2}\right\rfloor}\\
\vdots&\vdots&\vdots&&\vdots&\vdots\\
Q_{0}&\Im\left(\zeta^{N-2}\right)Q_{1}&\Re\left(\zeta^{N-2}\right)Q_{1}&\ldots&\Im\left(\zeta^{(N-2)\left\lfloor {N}/{2}\right\rfloor}\right)Q_{\left\lfloor {N}/{2}\right\rfloor}&\Re\left(\zeta^{(N-2)\left\lfloor {N}/{2}\right\rfloor}\right)Q_{\left\lfloor {N}/{2}\right\rfloor}\\
Q_{0}&\Im\left(\zeta^{N-1}\right)Q_{1}&\Re\left(\zeta^{N-1}\right)Q_{1}&\ldots&\Im\left(\zeta^{(N-1)\left\lfloor {N}/{2}\right\rfloor}\right)Q_{\left\lfloor {N}/{2}\right\rfloor}&\Re\left(\zeta^{(N-1)\left\lfloor {N}/{2}\right\rfloor}\right)Q_{\left\lfloor {N}/{2}\right\rfloor}
\end{array}\right).$$

Recall that in the configuration and momentum coordinates, the higher order Hamiltonian function is given by
$$H_2(q_i,p_i)=\sum_{i=1}^N \mu \left(q_{i1}^4+ q_{i2}^4\right).$$ Thus, we only need to investigate the configuration coordinates. Based on the structure of each $Q_j$, as shown in equations~\eqref{eq:qJnotEqualNover2} and~\eqref{eq:qJEqualsNover2}, they can be written as
\begin{equation}\label{eq:qi1Odd}
\begin{split}
	q_{i1}={}&QP[4(i-1)+1,:]\cdot X\\
	={}&Q_0[1,1]x_{01}+Q_0[1,4]x_{04}\\
	{}&+\sum_{j=1}^{\left\lfloor {N}/{2}\right\rfloor-1}\Im\left(\zeta^{(i-1)j}\right)\left(Q_j[1,1]x_{j1}+Q_j[1,4]x_{j4}\right)\\
	{}&+\Re\left(\zeta^{(i-1)j}\right)\left(Q_j[1,1]y_{j1}+Q_j[1,4]y_{j4}\right)\\
	{}&+\Im\left(\zeta^{(i-1)\left\lfloor {N}/{2}\right\rfloor}\right)\left(Q_{\left\lfloor {N}/{2}\right\rfloor}[1,3]x_{\left\lfloor {N}/{2}\right\rfloor 3}+Q_{\left\lfloor {N}/{2}\right\rfloor}[1,4]x_{\left\lfloor {N}/{2}\right\rfloor 4}\right)\\
	{}&+\Re\left(\zeta^{(i-1)\left\lfloor {N}/{2}\right\rfloor}\right)\left(Q_{\left\lfloor {N}/{2}\right\rfloor}[1,3]y_{\left\lfloor {N}/{2}\right\rfloor 3}+Q_{\left\lfloor {N}/{2}\right\rfloor}[1,4]y_{\left\lfloor {N}/{2}\right\rfloor 4}\right)
\end{split}
\end{equation}
and
\begin{equation}\label{eq:qi2Odd}
\begin{split}
	q_{i2}={}&QP[4(i-1)+2,:]\cdot X\\
	={}&Q_0[2,2]x_{02}+Q_0[2,3]x_{03}\\
	{}&+\sum_{j=1}^{\left\lfloor {N}/{2}\right\rfloor-1}\Im\left(\zeta^{(i-1)j}\right)\left(Q_j[2,2]x_{j2}+Q_j[2,3]x_{j3}\right)\\
	{}&+\Re\left(\zeta^{(i-1)j}\right)\left(Q_j[2,2]y_{j2}+Q_j[2,3]y_{j3}\right)\\
	{}&+\Im\left(\zeta^{(i-1)\left\lfloor N/2\right\rfloor}\right)\left(Q_{\left\lfloor {N}/{2}\right\rfloor}[2,1]x_{\left\lfloor {N}/{2}\right\rfloor 1}+Q_{\left\lfloor {N}/{2}\right\rfloor}[2,2]x_{\left\lfloor {N}/{2}\right\rfloor 2}\right)\\
	{}&+\Re\left(\zeta^{(i-1)\left\lfloor N/2\right\rfloor}\right)\left(Q_{\left\lfloor {N}/{2}\right\rfloor}[2,1]y_{\left\lfloor {N}/{2}\right\rfloor 1}+Q_{\left\lfloor {N}/{2}\right\rfloor}[2,2]y_{\left\lfloor {N}/{2}\right\rfloor 2}\right)
\end{split}
\end{equation}
where $QP[i,:]$ denotes the $i^{th}$ row and $QP[i,j]$ denotes the entry of the $i^{th}$ row and $j^{th}$ column  of the matrix $QP$.

Recall that the  $H_2(q_i,p_i)=\sum_{i=1}^N\frac{1}{2}\left(q_{i1}^4+q_{i2}^4\right)$ term represents the nonlinearities of the system in the Hamiltonian. From~\eqref{eq:qi1Odd} and~\eqref{eq:qi2Odd}, we see that $x_{01}$ only appears in $q_{i1}$ and $x_{03}$ only appears in $q_{i2}$. Products of $x_{01}$ and $x_{03}$ do not appear in $H_2$ nor do these two variables multiply any common factors. Since any degree four invariant involving $x_{01}$ and $x_{03}$ must contain their products, we can conclude that degree four invariants with these two variables do not appear in $H_2(X)$. Using the same reasoning, we can deduce that there are no degree four invariants with $x_{02}$ and $x_{04}$ terms as well. Thus, there are no degree four invariants involving $\mathcal U_1$ from~\eqref{eq:degree2InvariantsOdd}. By similar reasoning, we can rule out any degree four invariants involving $\mathcal U_{3}$ as well.

Based on these observations, degree four invariants must take the form of $\mathcal U_{2}^2$, but not all possible combinations of $\mathcal U_{2}^2$ are realized in $H_2$.  To simplify the notation, we divide the possible forms of $\mathcal U_2$ as
\begin{equation*}
\begin{split}
\mathcal U_{21}={}&x^2_{ja}+y^2_{ja},\\
\mathcal U_{22}={}&x^2_{\lfloor N/2\rfloor s}+y^2_{\lfloor N/2\rfloor s},\\
\mathcal U_{23}={}&x^2_{jb}+y^2_{jb},\AND\\
\mathcal U_{24}={}&x^2_{\lfloor N/2\rfloor t}+y^2_{\lfloor N/2\rfloor t},\\
\end{split}
\end{equation*} for $a=1,4$, $b=2,3$, $s=3,4$, $t=1,2$, and $j=1,\ldots, \lfloor N/2\rfloor-1$. We have divided the possible $\mathcal U_{2}$ terms in this manner because $\mathcal U_{21}$ and $\mathcal U_{22}$ correspond to the possible terms that can arise in $q_{i1}^4$. Similarly, $\mathcal U_{23}$ and $\mathcal U_{24}$ correspond to the possible terms that can occur in $q_{i2}^4$. Thus, when $N$ is odd,  the possible degree four invariants have the form
\begin{equation*}
\mathcal U_{21}^2,\; \mathcal U_{21}\mathcal U_{22},\; \mathcal U_{22}^2,\; \mathcal U_{23}^2,\; \mathcal U_{23}\mathcal U_{24},\;\mbox{and }\; \mathcal U_{24}^2.
\end{equation*}
In this notation, the multiplication of the $\mathcal U$ terms are done over all possible permutation of the indexes. For example, $(x^2_{11}+y^2_{11})(x^2_{34}+y^2_{34})$ is a possible product from $\mathcal U_{21}^2$ for $a=1,4$.

When $N$ is even, the $P$ may be written in block matrix form as
$$\footnotesize
P=\left(\begin{array}{ccccccccc}
I_4&I_4&I_4&\ldots&I_4&I_4&I_4\\
I_4&\Im(\zeta)I_4&\Re\left(\zeta\right)I_4&\ldots&\Im\left(\zeta^{\left({N}/{2}-1\right)}\right)I_4&\Re\left(\zeta^{\left({N}/{2}-1\right)}\right)I_4&-I_4\\
\vdots&\vdots&\vdots&&\vdots&\vdots\\
I_4&\Im\left(\zeta^{\left(N-1\right)}\right)I_4&\Re\left(\zeta^{\left(N-1\right)}\right)I_4&\ldots&\Im\left(\zeta^{\left(N-1\right)\left({N}/{2}-1\right)}\right)I_4&\Re\left(\zeta^{(N-1)\left({N}/{2}-1\right)}\right)I_4&-I_4
\end{array}\right)$$
and $Q$ may also be written as $Q=\operatorname {diag}\left(Q_0,Q_1,Q_1,\ldots,Q_{{N}/{2}-1},Q_{\frac{ N}{2}},Q_{{N}/{2}}\right)$. In this case, the product of $Q$ and $P$ is
$$\tiny
QP=\left(\begin{array}{ccccccccc}
Q_0&Q_1&Q_1&\ldots&Q_{{N}/{2}-1}&Q_{{N}/{2}-1}&Q_{{N}/{2}}\\
Q_0&\Im(\zeta)Q_1&\Re\left(\zeta\right)Q_1&\ldots&\Im\left(\zeta^{\left({N}/{2}-1\right)}\right)Q_{{N}/{2}-1}&\Re\left(\zeta^{\left({N}/{2}-1\right)}\right)Q_{{N}/{2}-1}&-Q_{{N}/{2}}\\
Q_0&\Im\left(\zeta^2\right)Q_1&\Re\left(\zeta^2\right)Q_1&\ldots&\Im\left(\zeta^{2\left({N}/{2}-1\right)}\right)Q_{{N}/{2}-1}&\Re\left(\zeta^{2\left({N}/{2}-1\right)}\right)Q_{{N}/{2}-1}&Q_{{N}/{2}}\\
\vdots&\vdots&\vdots&&\vdots&\vdots&\vdots\\
Q_{0}&\Im\left(\zeta^{(N-2)}\right)Q_{1}&\Re\left(\zeta^{(N-2)}\right)Q_{1}&\ldots&\Im\left(\zeta^{(N-2)\left({N}/{2}-1\right)}\right)Q_{{N}/{2}-1}&\Re\left(\zeta^{(N-2)\left({N}/{2}-1\right)}\right)Q_{{N}/{2}-1}&Q_{{N}/{2}}\\
Q_{0}&\Im\left(\zeta^{\left(N-1\right)}\right)Q_{1}&\Re\left(\zeta^{\left(N-1\right)}\right)Q_{1}&\ldots&\Im\left(\zeta^{(N-1)\left({N}/{2}-1\right)}\right)Q_{{N}/{2}-1}&\Re\left(\zeta^{(N-1)\left({N}/{2}-1\right)}\right)Q_{{N}/{2}-1}&-Q_{{N}/{2}}
\end{array}\right).$$

As in the case when $N$ is odd, we only need to investigate the configuration coordinates because of the form of $H_2$. Based on the structure of each $Q_j$, as shown in equations~\eqref{eq:qJnotEqualNover2} and~\eqref{eq:qJEqualsNover2}, they can be written as
\begin{align*}
	q_{i1}={}&QP[4(i-1)+1,:]\cdot X\\
	={}&Q_0[1,1]x_{01}+Q_0[1,4]x_{04}\\
	{}&+\sum_{j=1}^{ {N}/{2}-1}\Im\left(\zeta^{(i-1)j}\right)\left(Q_j[1,1]x_{j1}+Q_j[1,4]x_{j4}\right)\\
	{}&+\Re\left(\zeta^{(i-1)j}\right)\left(Q_j[1,1]y_{j1}+Q_j[1,4]y_{j4}\right)\\
	{}&+(-1)^{i-1}\left(Q_{ {N}/{2}}[1,3]x_{(N/2)3}+Q_{ {N}/{2}}[1,4]x_{(N/2) 4}\right)
\end{align*}
and
\begin{align*}
	q_{i2}={}&QP[4(i-1)+2,:]\cdot X\\
	={}&Q_0[2,2]x_{02}+Q_0[2,3]x_{03}\\
	{}&+\sum_{j=1}^{ {N}/{2}-1}\Im\left(\zeta^{(i-1)j}\right)\left(Q_j[2,2]x_{j2}+Q_j[2,3]x_{j3}\right)\\
	{}&+\Re\left(\zeta^{(i-1)j}\right)\left(Q_j[2,2]y_{j2}+Q_j[2,3]y_{j3}\right)\\
	{}&+(-1)^{i-1}\left(Q_{{N}/{2}}[2,1]x_{(N/2) 1}+Q_{ {N}/{2}}[2,2]x_{(N/2)2}\right)
\end{align*}
where $QP[i,:]$ denotes the $i^{th}$ row and $QP[i,j]$ denotes the entry of the $i^{th}$ row and $j^{th}$ column  of the matrix $QP$.

For reasons stated in the case when $N$ is odd, degree four invariants involving $\mathcal U_1$ do not occur. Again, we observe that not all possible forms of $\mathcal U_2$ and $\mathcal U_4$ can be realized. We divide the possible $\mathcal U_2$ and $\mathcal U_4$ as follow:
\begin{equation*}
\begin{split}
\mathcal U_{21}={}&x^2_{ja}+y^2_{ja},\\
\mathcal U_{41}={}&x^2_{(N/2)s},\\
\mathcal U_{42}={}&x_{(N/2)3}x_{(N/2)4},\\
\mathcal U_{23}={}&x^2_{jb}+y^2_{jb},\\
\mathcal U_{43}={}&x^2_{(N/2)t},\AND\\
\mathcal U_{44}={}&x_{(N/2)1}x_{(N/2)2},\\
\end{split}
\end{equation*} where $a=1,4$, $b=2,3$, $s=3,4$, $t=1,2$, and $j=1,\ldots, N/2-1$. The possible terms are divided because $\mathcal U_{21}$, $\mathcal U_{41}$ and $\mathcal U_{42}$ only appear in $q_{I1}^4$ and $\mathcal U_{23}$, $\mathcal U_{43}$ and $\mathcal U_{44}$ only appear in $q_{I2}^4$. Thus, when $N$ is even, the possible degree four invariants are
\begin{equation*}
\mathcal U_{21}^2, \;\mathcal U_{21}\mathcal U_{41},\;\mathcal U_{21}\mathcal U_{42},\; \mathcal U_{41}^2,\; \mathcal U_{42}^2, \;\mathcal U_{41}\mathcal U_{42}, \;\mathcal U_{23}^2,\;\mathcal U_{23}\mathcal U_{43},\;\mathcal U_{23}\mathcal U_{44}, \;\mathcal U_{43}^2,\;\mathcal U_{44}^2,\;\mbox{and}\;\;\mathcal U_{43}\mathcal U_{44}.
\end{equation*}
These products are multiplied over all possible combinations of the indexes as we noted in the case when $N$ is odd.

\newpage
\section*{References}
\bibliography{gyroscope}
\bibliographystyle{iopart-num}
\end{document}